\documentclass[a4paper, 10pt, twoside, notitlepage]{amsart}

\usepackage{amsmath,amscd}
\usepackage{amssymb}
\usepackage{amsthm}
\usepackage{comment}
\usepackage{graphicx, xcolor}

\usepackage{mathrsfs}
\usepackage[ocgcolorlinks, linkcolor=blue]{hyperref}

\usepackage{bm}
\usepackage{bbm}
\usepackage{url}

\usepackage[utf8]{inputenc}
\usepackage{mathtools,amssymb}
\usepackage{esint}
\usepackage{tikz}
\usepackage{dsfont}
\usepackage{relsize}
\usepackage{url}
\usepackage{xcolor}
\usepackage{graphicx}
\usepackage{mathrsfs}
\usepackage[shortlabels]{enumitem}
\usepackage{lineno}
\usepackage{amsmath}
\usepackage{enumitem}
\usepackage{amsthm} 
\usepackage{verbatim}
\usepackage{dsfont}
\numberwithin{equation}{section}

\allowdisplaybreaks

\mathtoolsset{showonlyrefs}

\graphicspath{{images/}}

\newtheorem{theorem}{Theorem}[section]
\newtheorem{lemma}[theorem]{Lemma}
\newtheorem{definition}[theorem]{Definition}
\newtheorem{corollary}[theorem]{Corollary}

\newtheorem{claim}[theorem]{Claim}
\newtheorem{proposition}[theorem]{Proposition}
\newtheorem{example}[theorem]{Example}
\newtheorem{problem}[theorem]{Problem}
\newtheorem{question}{Question}
\newtheorem{remark}[theorem]{Remark}
\newtheorem{assumption}[theorem]{Assumption}

\title[Calder\'on problem for nonlocal viscous wave equations]{Calder\'on problem for nonlocal viscous wave equations: Unique determination of linear and nonlinear perturbations}

\author[P. Zimmermann]{Philipp Zimmermann}
\address{Department of Mathematics, ETH Zurich, Z\"urich, Switzerland \& Departament de Matem\`atiques i Inform\`atica, Universitat de Barcelona, Barcelona, Spain}
\email{philipp.zimmermann@ub.edu}

\newcommand{\C}{{\mathbb C}}
\newcommand{\R}{{\mathbb R}}

\newcommand{\N}{{\mathbb N}}

\newcommand{\eps}{\varepsilon}

\newcommand {\p} {\partial}

\newcommand{\LC}{\left(}
\newcommand{\RC}{\right)}

\newcommand{\schwartz}{\mathscr{S}}

\newcommand{\tempered}{\mathscr{S}^{\prime}}

\newcommand{\fourier}{\mathcal{F}}
\newcommand{\ifourier}{\mathcal{F}^{-1}}

\newcommand{\distr}{\mathscr{D}^{\prime}}

\newcommand{\abs}[1]{\left\lvert #1 \right\rvert}
\DeclareMathOperator{\Div}{div} 
\DeclareMathOperator{\supp}{supp} 
\DeclareMathOperator{\loc}{loc} 

\newcommand{\weak}{\rightharpoonup}

\newcommand{\Vareps}{\boldsymbol{\varepsilon}}

\begin{document}
	
	\maketitle
	\begin{abstract}
		The main goal of this article is to study a Calder\'on type inverse problem for certain viscous nonlocal wave equations. We show that the partial Dirichlet to Neumann map uniquely determines on the one hand linear perturbations and on the other hand homogeneous nonlinearities $f(u)$ whenever the latter satisfy a certain growth assumption. As a preliminary step we discuss the well-posedness in each case, where for the nonlinear setting we invoke the implicit function theorem after establishing the differentiability of the associated Nemytskii operator $f(u)$. In the linear case we establish a Runge approximation theorem in $L^2(0,T;\widetilde{H}^{s}(\Omega))$, which allows us to uniquely determine potentials that belong only to $L^{\infty}(0,T;L^p(\Omega))$ for some $1<p\leq \infty$ satisfying suitable restrictions. In the nonlinear case, we first derive an appropriate integral identity and combine this with the differentiability of the solution map around zero to show that the nonlinearity is uniquely determined by the Dirichlet to Neumann map. To make this linearization technique work, it is essential that we have a Runge approximation in $L^2(0,T;\widetilde{H}^s(\Omega))$ instead of $L^2(\Omega_T)$ at our disposal.
		
		\medskip
		
		

	\end{abstract}

	\tableofcontents

	\section{Introduction}
	\label{sec: introduction}

	In recent years, many different inverse problems for nonlocal partial differential equations (PDEs) have been studied in the literature. The very first work in this area was the article \cite{GSU20} by Ghosh, Salo and Uhlmann. They showed that the nonnegative potential $q\in L^{\infty}(\Omega)$ in the \emph{fractional Schr\"odinger equation}
	\begin{equation}
		\label{eq: fractional Schroedinger equation}
		((-\Delta)^s+q)u=0\text{ in }\Omega
	\end{equation}
	is uniquely determined by the \emph{(partial) Dirichlet to Neumann (DN) map}
	\begin{equation}
		\Lambda_q f= (-\Delta)^s u_f|_{W_2}, \quad f\in C_c^{\infty}(W_1),
	\end{equation}
	for arbitrary fixed measurement sets $W_1,W_2\subset \Omega_e$. Here, $\Omega\subset\R^n$ is a bounded domain with exterior $\Omega_e=\R^n\setminus\overline{\Omega}$, $0<s<1$ and $u_f\colon\R^n\to\R$ denotes the unique solution to \eqref{eq: fractional Schroedinger equation} with exterior data $u_f|_{\Omega_e}=f$. An essential analytical tool in their work is the so-called \emph{unique continuation property (UCP)} of the fractional Laplacian. Roughly speaking, the UCP can be phrased as follows:
	
	\medskip
	
	\begin{center}
		If $u:\R^n\to \R$ satisfies $(-\Delta)^su=u=0$ in an open set $V\subset \R^n$, then $u\equiv 0$.\vspace{0.1cm}\\
	\end{center}
	
	\noindent
	Its proof depends on the famous \emph{Caffarelli--Silvestre (CS) extension} \cite{CS-extension-problem-fractional-laplacian} of the fractional Laplacian, which allows to characterize the fractional Laplacian of $u$ as the Neumann data of the solution $U$ to the degenerate elliptic PDE
	\begin{equation}
		\label{eq: CS extension PDE}
		\Div \LC y^{1-2s}\nabla U(x,y)\RC=0 \text{ in }\R^{n+1}_+
	\end{equation}
	with Dirichlet data $U(x,0)=u$ on $\partial \R^{n+1}_+$. Solutions to such equations, having $A_2$ Muckenhoupt weights as coefficients, have already been studied a long time ago in the celebrated work \cite{fabes1982local} by Fabes, Kenig and Serapioni. Let us note that CS type extensions are only available for a restricted classes of nonlocal operators as discussed in more detail by Kwa\'snicki, Mucha and Stinga, Torrea in \cite{kwasnicki2018extension} and \cite{ST10}, respectively. Based on this fact, in subsequent research articles in this field, the main focus was put on nonlocal inverse problems for equations of the form
	\begin{equation}
		\label{eq: generalized frac cal}
		L_K u+Q(u)=0\text{ in }\Omega,
	\end{equation}
	where $L_K$ is an elliptic (with potentially variable coefficients) nonlocal operator having the UCP and $Q$ is a possibly nonlinear function of $u$. Prototypical results in this field show that the DN map associated with equation \eqref{eq: generalized frac cal} uniquely determines the function $Q$ (see e.g.~\cite{bhattacharyya2021inverse,cekic2020calderon,harrach2017nonlocal-monotonicity,GRSU18,GU2021calder,RS17,RZ-unbounded}). Additionally, there are works in which the authors have attempted to recover the coefficients $K$, on which the nonlocal operator $L_K$ may depend, from the associated DN map. Examples of such nonlocal operators $L$ include the fractional powers of elliptic second order operators $\mathrm{L}^s_{\sigma}=\LC -\Div(\sigma\nabla)\RC ^s$ and the fractional conductivity operator $L^s_{\gamma}$. Both examples fall into the class of elliptic integro-differential operators, denoted by $\mathcal{L}_0$, which consists of all operators $L_K$ that can be written in strong form as
    \begin{equation}
        L_Ku(x)=\text{p.v.}\int_{\R^n}K(x,y)(u(x)-u(y))\,dy,
    \end{equation}
    where the kernel $K(x,y)$ satisfies
    \begin{equation}
    \label{eq: elliptic integro diff}
        K(x,y)=K(y,x)\quad\text{and}\quad \frac{c}{|x-y|^{n+2s}}\leq K(x,y)\leq \frac{C}{|x-y|^{n+2s}}.
    \end{equation}
    The second condition in \eqref{eq: elliptic integro diff} simply means that the kernel $K(x,y)$ is comparable to that of the fractional Laplacian $(-\Delta)^s$
    (see Section~\ref{sec: preliminaries} for more details on the fractional Laplacian). In \cite[Section 2.3]{GLX} it is shown that for every uniformly elliptic coefficient $\sigma\in C^{\infty}(\R^n)$, the nonlocal operator $\mathrm{L}_{\sigma}^s=\LC -\Div(\sigma\nabla)\RC ^s$ belongs to $\mathcal{L}_0$, and its kernel is given by
    \[
        \mathrm{K}^s_{\sigma}(x,y)=\frac{1}{\Gamma(-s)}\int_0^{\infty}p_t(x,y)\frac{dt}{t^{1+s}},
    \]
    where $\Gamma$ denotes Euler's Gamma function and $p_t$ is the symmetric heat kernel associated with the elliptic operator $\mathrm{L}_{\sigma}=-\Div(\sigma\nabla)$. We note that the bounds in \eqref{eq: elliptic integro diff} for $\mathrm{K}^s_{\sigma}$ follow from the fact that the heat kernel $p_t$ can be controlled from both sides by Gaussians; that is,
    \[
        c (4\pi t)^{-n/2} e^{-\frac{|x-y|^2}{4\alpha t}}\leq  p_t(x,y)\leq C (4\pi t)^{-n/2} e^{\frac{|x-y|^2}{4\beta t}}
    \]
    for some constants $c,C,\alpha,\beta>0$.
    To see that $L^s_{\gamma}\in\mathcal{L}_0$, we recall that the kernel of the fractional conductivity operator $L^s_{\gamma}$ is given by
    \[
        K^s_{\gamma}(x,y)=C_{n,s}\frac{\gamma^{1/2}(x)\gamma^{1/2}(y)}{|x-y|^{n+2s}}
    \] 
    where $C_{n,s}>0$ is a suitable normalization constant and $\gamma\colon\R^n\to\R$ is a uniformly elliptic function. The works \cite{RZ2022LowReg, ruland2023revisiting} provide affirmative answers to the question of whether the kernels $K^s_{\gamma}$ and $\mathrm{K}^s_{\sigma}$ can be recovered from the associated DN maps. Let us remark that there are only very few results on the simultaneous recovery of the leading order coefficient $K$ and the lower order perturbation $Q$. For this type of questions, we may point the interested reader to the works \cite{zimmermann2023inverse}, \cite{LZ2024uniqueness}, and \cite{Trans-anisotropic-LNZ}.
    
    Later on these studies were extended to the parabolic setting (e.g., \cite{LRZ2022calder,LLR2019calder,LLU2022para,LLU2023calder,LZ2023unique}). Recently, Kow, Lin and Wang studied in \cite{KLW2022} a Calder\'on type inverse problem related to the \emph{nonlocal wave equation}\footnote{For any set $A\subset\R^n$ and $T>0$, we write $A_T$ to denote the space time cylinder $A\times (0,T)$.}
    \begin{equation}
    \label{eq: nonlocal wave eq}
        \begin{cases}
				\LC \partial_t^2 +(-\Delta)^s+q\RC u=0 &\text{ in }\Omega_T\\
				u=\varphi     &\text{ in }(\Omega_e)_T,\\
				u(0)=0,\quad \partial_t u(0)=0  &\text{ in }\Omega
			\end{cases}
    \end{equation}
    with $0<s<1$, $q=q(x)\in L^{\infty}(\Omega)$ and showed under suitable assumptions on the domains that:
    \begin{enumerate}[(i)]
        \item\label{item: Runge wave eq} \textit{Runge approximation:} Any $v\in L^2(\Omega_T)$ can be approximated arbitrarily well in $L^2(\Omega_T)$ by solutions $u_{\varphi}|_{\Omega_T}$ of \eqref{eq: nonlocal wave eq}, where $\varphi\in C_c^{\infty}(W_T)$ for some fixed measurement set $W\Subset \Omega_e$.
        \item\label{item: unique det wave eq} \textit{Unique determination:} Let $q_1,q_2\in L^{\infty}(\Omega)$ be two potentials and denote by $\Lambda_{q_j}$ its DN map related to \eqref{eq: nonlocal wave eq}, i.e.
        \[
            \Lambda_{q_j}\varphi=\left.(-\Delta)^su_{\varphi}\right|_{(\Omega_e)_T}.
        \] 
        If one has 
        \[
            \left.\Lambda_{q_1}\varphi\right|_{(W_2)_T}=\left.\Lambda_{q_2}\varphi\right|_{(W_2)_T}
        \]
        for all $\varphi\in C_c^{\infty}((W_1)_T)$, where $W_1,W_2\subset\Omega_e$ are two fixed measurement sets, then one has
        \[
            q_1(x)=q_2(x)\text{ in }\Omega.
        \]
    \end{enumerate}
    Let us point out that in contrast to the elliptic (see \cite{RZ-unbounded}) or the parabolic case (see \cite{LRZ2022calder}), it is not known whether the Runge approximation for the nonlocal wave equation  \eqref{eq: nonlocal wave eq} holds in $L^2(0,T;\widetilde{H}^s(\Omega))$\footnote{The precise definition of this space is given in Section~\ref{sec: preliminaries}.}.Hence, the following question remains open.
    \begin{question}
    \label{Runge quest}
        Is the Runge set
    \[
        \mathcal{R}_W=\{u_{\varphi}-\varphi\,;\,\varphi\in C_c^{\infty}(W_T)\},
    \]
    where the notation $u_{\varphi}$ is as above, dense in $L^2(0,T;\widetilde{H}^s(\Omega))$?
    \end{question}
    The main obstruction in proving such a Runge approximation is the low regularity of solutions to the equation 
    \[
    \begin{cases}
				\LC \partial_t^2 +(-\Delta)^s+q\RC u=F &\text{ in }\Omega_T\\
				u=0     &\text{ in }(\Omega_e)_T,\\
				u(0)=0,\quad \partial_t u(0)=0 &\text{ in }\Omega,
			\end{cases}
   \]
   when $F$ only belongs to the space $L^2(0,T;H^{-s}(\Omega))$ instead of $L^2(\Omega_T)$. This phenomenon already appears in the local case $s=1$ (see \cite[Chapter 11]{Precup+2013}). On the other hand, a main difference between the local and nonlocal wave equation \eqref{eq: nonlocal wave eq} is that in the latter case the equation does not have a finite speed of propagation, which in turn relies on the UCP of the fractional Laplacian. 

    Because of the lack of such a density result the techniques of this article cannot be directly adapted to the nonlinear nonlocal wave equation and it is an open question, whether the DN map related to the \emph{nonlinear nonlocal wave equation}
     \begin{equation}
    \label{eq: nonlinear nonlocal wave eq}
        \begin{cases}
				\LC \partial_t^2 +(-\Delta)^s\RC u+f(u)=0 &\text{ in }\Omega_T\\
				u=\varphi     &\text{ in }(\Omega_e)_T,\\
				u(0)=0,\quad \partial_t u(0)=0 &\text{ in }\Omega
			\end{cases}
    \end{equation}
    uniquely determines suitable nonlinearities $f$. Here and below, $f(u)$ denotes the Nemytskii operator associated to a Carath\'eodory function $f\colon \Omega\times\R\to \R$, that is 
   \begin{equation}
   \label{eq: Nemytskii op intro}
       f(u)(x,t)=f(x,u(x,t)).
   \end{equation}

    In this article, we study a Calder\'on type inverse problem for linear and nonlinear perturbations of the \emph{nonlocal viscous wave equation}
     \begin{equation}
    \label{eq: nonlocal viscous wave eq}
        \begin{cases}
				\LC \partial_t^2 +(-\Delta)^s\partial_t+(-\Delta)^s\RC u=0 &\text{ in }\Omega_T\\
				u=\varphi     &\text{ in }(\Omega_e)_T,\\
				u(0)=0,\quad \partial_t u(0)=0 &\text{ in }\Omega.
			\end{cases}
    \end{equation}
    The terminology used for this equation is discussed in the next section. More concretely, this means that we study the problems
    \begin{equation}
    \label{eq: nonlocal viscous wave eq lin pertur}
        \begin{cases}
				\LC \partial_t^2 +(-\Delta)^s\partial_t+(-\Delta)^s+q\RC u=0 &\text{ in }\Omega_T\\
				u=\varphi     &\text{ in }(\Omega_e)_T,\\
				u(0)=0,\quad \partial_t u(0)=0 &\text{ in }\Omega
			\end{cases}
    \end{equation}
    and
    \begin{equation}
    \label{eq: nonlocal viscous wave eq nonlin pertur}
        \begin{cases}
				\LC \partial_t^2 +(-\Delta)^s\partial_t+(-\Delta)^s\RC u+f(u)=0 &\text{ in }\Omega_T\\
				u=\varphi     &\text{ in }(\Omega_e)_T,\\
				u(0)=0,\quad \partial_t u(0)=0 &\text{ in }\Omega
			\end{cases}
    \end{equation}
    and aim to uniquely recover the potential $q$ and nonlinearity $f$ under suitable assumptions from the related DN maps
   \begin{equation}
   \label{eq: DN map lin and nonlin}
   \begin{split}
       \Lambda_q \varphi &=\left.((-\Delta)^su_\varphi+(-\Delta)^s\partial_t u_{\varphi})\right|_{(\Omega_e)_T}\\
       \Lambda_{f} \varphi &=\left.((-\Delta)^sv_\varphi+(-\Delta)^s\partial_t v_{\varphi})\right|_{(\Omega_e)_T},
   \end{split}
   \end{equation}
   where $u_\varphi, v_\varphi$ are the unique solutions to \eqref{eq: nonlocal viscous wave eq lin pertur} and \eqref{eq: nonlocal viscous wave eq nonlin pertur}, respectively. Let us remark that in contrast to the results in \cite{KLW2022} the potential in \eqref{eq: nonlocal viscous wave eq lin pertur} is allowed to vary in time and is not necessarily bounded in space.

    Finally, let us note that in fact one can construct unique solutions to the linear nonlocal wave equations by first considering solutions $u_\eps$ to the nonlocal viscous wave equation with loss term $\eps(-\Delta)^s\partial_t$ and then passing to the limit $\eps\to 0$ (see \cite[Chapter XVIII]{DautrayLionsVol5}).

   \subsection{Nonlocal viscous wave equations and related models}

    The main goal of this section is to motivate the terminology for equation \eqref{eq: nonlocal viscous wave eq} and to discuss related models. 

    First of all let us recall that the initial and boundary value problem for the \emph{viscous wave equation} is given by
    \begin{equation}
    \label{eq: viscous wave eq}
        \begin{cases}
				\LC \frac{1}{c^2}\partial_t^2 +\tau(-\Delta)\partial_t+(-\Delta)\RC u=0 &\text{ in }\Omega_T\\
				u=\varphi     &\text{ on }\partial\Omega_T,\\
				u(0)=0,\quad \partial_t u(0)=0 &\text{ in }\Omega,
			\end{cases}
    \end{equation}
    which emerges in acoustics to describe the propagation of sound in a viscous fluid. The quantity $u$ represents the sound pressure, $c$ the speed of propagation and $\tau$ the relaxation time, which can be calculated as 
    \[
        \tau = \frac{4\mu}{3\rho_0c^2}
    \]
    with $\mu$ being the shear bulk viscosity coefficient, which has been measured for many fluids, and $\rho_0$ is the static density. The term $ \tau(-\Delta)u$ corresponds physically to an additional loss term. If we formally put $c=\tau=1$ and replace the Laplacian by the fractional Laplacian $(-\Delta)^s$, then we arrive at the problem \eqref{eq: nonlocal viscous wave eq}.

    Next, we describe a time fractional generalization of the problem \eqref{eq: viscous wave eq} and a generalization of \eqref{eq: nonlocal viscous wave eq}.
    
   \begin{enumerate}[{(G}1)]
       \item\label{model 1} The first model we introduce reads as follows:
       \begin{equation}
       \label{eq: second generalization}
           \begin{cases}
				\LC \frac{1}{c^2}\partial_t^2 +\beta\tau^\beta(-\Delta)\partial_t^\beta+(-\Delta)\RC u=0 &\text{ in }\Omega_T\\
				u=\varphi     &\text{ on }\partial\Omega_T,\\
				u(0)=0,\quad \partial_t u(0)=0 &\text{ in }\Omega.
			\end{cases}
       \end{equation}
      Here \( \partial_t^{\beta} \) denotes a fractional time derivative. Such operators naturally arise when interpolating between Hooke's law—which relates the strain and stress in an elastic solid—and Newton's law for fluids, which describes the linear relationship between stress and the strain rate in an ideal viscous fluid.  For further information on this model, we refer to \cite{wang2016generalized, xu2023determination} and the references therein.
       \item\label{model 2} An immediate generalization of the model \eqref{eq: nonlocal viscous wave eq} is
       \begin{equation}
       \label{eq: first generalization}
           \begin{cases}
				\LC \partial_t^2 +(-\Delta)^{s_1}\partial_t+(-\Delta)^{s_2}\RC u=0 &\text{ in }\Omega_T\\
				u=\varphi     &\text{ in }(\Omega_e)_T,\\
				u(0)=0,\quad \partial_t u(0)=0 &\text{ in }\Omega,
			\end{cases}
       \end{equation}
       where two fractional Laplacians of different orders \( s_1 \) and \( s_2 \) appear. 
    The special case \( s_1 = \tfrac{1}{2} \) and \( s_2 = 1 \) has recently attracted interest and has been studied in the setting \( \Omega = \mathbb{R}^2 \), for example in \cite{kuan2021deterministic, kuan2022probabilistic, deroubin2023norm}. 
    In these works, the equation may also include a nonlinearity of the form \( |u|^{p-1}u \) for \( p>1 \) or \( u^k \) for \( k \ge 2 \). 
    Such power-type nonlinearities have been extensively studied in the context of dispersive wave equations, such as the nonlinear wave equation or the nonlinear Schr\"odinger equation
       \[
            i\partial_t u+\Delta u+f(u)=0\text{ in }\R^n
       \]
       (see \cite{tao2006nonlinear}).  Inverse problems related to the model \eqref{eq: first generalization} are addressed in upcoming work of the author with Katya Krupchyk \cite{PZ-KK}.
   \end{enumerate}
    
    \subsection{Main results}
    
   In this section, we present our main results on the inverse problems for the nonlocal viscous wave equation with linear and nonlinear perturbations.
    
	\begin{theorem}[Uniqueness of linear perturbations]
    \label{Thm: main linear}
		Let $\Omega\subset\R^n$ be a bounded Lipschitz domain, $T>0$ and $s>0$ a non-integer. Suppose that for $j=1,2$ we have given potentials $q_j\in L^1_{loc}(\Omega_T)$ such that
        \begin{enumerate}[(i)]
        \item $q_j\in L^{\infty}(0,T;L^p(\Omega))$ for some $1\leq p<\infty$ satisfying
        \[
		\begin{cases}
			n/s\leq p\leq \infty, &\, \text{if }\, 2s< n,\\
			2<p\leq \infty,  &\, \text{if }\, 2s= n,\\
			2\leq p\leq \infty, &\, \text{if }\, 2s > n,
		\end{cases}
		\]
            \item $t\mapsto \int_{\Omega}q_j(x,t)\varphi(x)\,dx\in C([0,T])$ for any $\varphi\in C_c^{\infty}(\Omega)$.
        \end{enumerate}
        Furthermore, assume that $W_1,W_2\subset\Omega_e$ are given measurement sets such that the DN maps $\Lambda_{q_j}$ related to 
        \[
        \begin{cases}
				\LC \partial_t^2 +(-\Delta)^s\partial_t +(-\Delta)^s+q_j\RC u=0 &\text{ in }\Omega_T\\
				u=\varphi     &\text{ in }(\Omega_e)_T,\\
				u(0)=0,\quad \partial_t u(0)=0 &\text{ in }\Omega.
			\end{cases}
        \]
        satisfy 
        \begin{equation}
        \label{eq: cond DN map linear}
            \Lambda_{q_1}\varphi=\Lambda_{q_2}\varphi\text{ in }(W_2)_T
        \end{equation}
        for all $\varphi\in C_c^{\infty}((W_1)_T)$. Then there holds 
        \begin{equation}
        \label{eq: equality of potentials}
            q_1(x,t)=q_2(x,t)\text{ in }\Omega_T.
        \end{equation}
	\end{theorem}

    In Section~\ref{subsec: unique determination linear}, we first present the proof of Theorem~\ref{Thm: main linear} for time-reversal invariant potentials, based on a suitable integral identity. Throughout this article, we define the \emph{time reversal} of a function \( Q \in L^1_{\mathrm{loc}}(V_T) \), where \( V \subset \mathbb{R}^n \) is an arbitrary open set, by
\begin{equation}
\label{eq: time reversal}
    Q^{\star}(x,t) = Q(x, T - t),
\end{equation}
and we say that \( Q \) is \emph{time-reversal invariant} if \( Q^\star = Q \).
We then establish the proof of Theorem~\ref{Thm: main linear} for general potentials.

Afterwards, in Section~\ref{IP nonlinear}, we extend this approach to show that nonlinear perturbations \( f(x,\tau) \) are uniquely determined by the DN map:

    \begin{theorem}[Uniqueness of nonlinear perturbations]
    \label{Thm: main nonlinear}
		Let $\Omega\subset\R^n$ be a bounded Lipschitz domain, $T>0$ and $s>0$ a non-integer. Suppose that for $j=1,2$ we have given nonlinearities $f_j\colon \Omega\times\R\to\R$ satisfying Assumption~\ref{main assumptions on nonlinearities} with $0<r\leq 2$ and $f_j$ is $r+1$ homogeneous in the second variable. Let $U^j_0\subset\widetilde{W}^s_{rest}((\Omega_e)_T),U^j_1\subset \widetilde{W}_{ext}(0,T;\widetilde{H}^s(\Omega))$ be the neighborhoods of the origin provided by Theorem~\ref{prop: differentiability of solution map}, which have the property that, for any $\varphi\in U^j_0$, the problem 
        \begin{equation}
        \label{PDE main thm nonlinear}
		\begin{cases}
			(\partial_t^2 +(-\Delta)^s\partial_t +(-\Delta)^s)u+f_j(u)= 0&\text{ in }\Omega_T\\
			u=\varphi&\text{ in }(\Omega_e)_T,\\
			u(0)=0,\quad \partial_t u(0)=0  &\text{ in }\Omega
		\end{cases}
		\end{equation}
        admits a unique solution $u\in U_1^j$.
        Furthermore, assume that $W_1,W_2\subset\Omega_e$ are given measurement sets such that the DN maps $\Lambda_{f_j}$ related to \eqref{PDE main thm nonlinear} satisfy 
        \begin{equation}
        \label{eq: cond DN map nonlinear}
            \Lambda_{f_1}\varphi=\Lambda_{f_2}\varphi\text{ in }(W_2)_T
        \end{equation}
        for all $\varphi\in U_0^1\cap U_0^2$ that are supported in $(W_1)_T$. Then there holds 
        \begin{equation}
        \label{eq: equality nonlinearities}
            f_1(x,\rho)=f_2(x,\rho)\text{ for }x\in\Omega\text{ and }\rho\in\R.
        \end{equation}
	\end{theorem}

    \begin{remark}
        Both uniqueness theorems can be extended to other nonlocal operators $L$ instead of the fractional Laplacian as long as they satisfy appropriate structural assumptions and a corresponding UCP. For this purpose we recall the necessary tools to solve the forward problem in a general framework in Section~\ref{subsec: abstract framework}. Note that this has been done in the case of elliptic nonlocal inverse problems in the work \cite{RZ-unbounded}.
    \end{remark}

 \subsection{Recent developments in inverse problems for nonlocal wave equations}

Since the appearance of the present article, the author, in collaboration with Y.-H.~Lin, S.-R.~Fu, T.~Tyni, and Y.~Yu, has obtained several unique determination results for the Calderón problem associated with nonlocal wave equations of the form
\begin{equation}
\label{eq: nonlocal wave}
    (\partial_t^2 + \gamma \partial_t + \lambda (-\Delta)^s \partial_t + (-\Delta)^s)u + f(u) = 0
\end{equation}
and third order nonlocal wave equations of the form
\begin{equation}
\label{eq: semilinear MGT eq}
    (\partial_t^3 + \alpha \partial_t^2 + b(-\Delta)^s \partial_t + c(-\Delta)^s)u + f(u) = 0.
\end{equation}

\begin{enumerate}[(i)]
    \item\label{gen 1}
    The works \cite{LTZ2} and \cite{FYZ-third-order} establish the recovery of potentials \( q \) and homogeneous nonlinearities \( f \) in the nonlocal wave equation \eqref{eq: nonlocal wave} with \( \gamma = \lambda = 0 \), as well as in the third order nonlocal wave equation \eqref{eq: semilinear MGT eq}, from the associated DN maps.  
    Their proofs rely on the classical \( L^2(\Omega_T) \) Runge approximation theorem for these equations.

    \item\label{gen 2}
    In \cite{LTZ1}, an \( L^2(0,T;\widetilde{H}^s(\Omega)) \) Runge approximation theory was developed for the nonlocal wave equation \eqref{eq: nonlocal wave} with \( \gamma = \lambda = 0 \).  
    This result was subsequently extended to the damped case \( \gamma \neq 0 \).  
    Consequently, these works resolve Question~\ref{Runge quest} raised on p.~3.  
    The key ingredient is the construction of so-called \emph{very weak solutions}, which allow rough initial data and sources.

    \item\label{gen 3}
    The improved Runge approximation results in \ref{gen 2} allow for the unique determination of polyhomogeneous nonlinearities \( f(x,u) \) whose growth at infinity exceeds the classes considered in earlier articles.
    Moreover, in \cite{PZ2} (for the equation \eqref{eq: nonlocal wave} with $\lambda=0$), the damping coefficient \( \gamma \) is permitted to vary over \( \Omega \), and it is shown that the DN map simultaneously determines both the nonlinearity \( f \) and the damping coefficient \( \gamma \).

    \item
    Finally, in \cite{FYZ-third-order} the nonlinearity \( f \) in \eqref{eq: semilinear MGT eq} is assumed to satisfy the same structural conditions as in the works above, but may additionally depend on the time variable \( t \).  
    The analysis developed there extends naturally to the other settings, and in particular Theorem~\ref{Thm: main nonlinear} can be generalized to nonlinearities of the form \( f(x,t,\tau) \).  
    For simplicity, however, we confine ourselves in this article to the time-independent case.
\end{enumerate}

	\subsection{Organization of the article}
    This article is organized as follows. In Section \ref{sec: preliminaries}, we introduce the functional analytic setup and in particular introduce the fractional Sobolev spaces, the fractional Laplacian and the Bochner Lebesgue spaces. Moreover, in the last paragraph of this section we discuss an abstract framework for solving some classes of second order in time PDEs. In Section ~\ref{sec: well-posedness}, we study the well-posedness theory of the problems \eqref{eq: nonlocal viscous wave eq lin pertur} and \eqref{eq: nonlocal viscous wave eq nonlin pertur}. The well-posedness of the nonlinear problem is achieved by invoking the implicit function theorem and to this end we study first in Section~\ref{subsec: diff of nemytskii} the differentiability of the Nemytskii operator $u\mapsto f(u)$ (Lemma~\ref{lemma: Differentiability of Nemytskii operators}). Then in Section~\ref{sec: solving linear IP}, after establishing the Runge approximation (Proposition~\ref{prop: runge}) and a suitable integral identity (Lemma~\ref{lemma: integral identity}), we prove Theorem~\ref{Thm: main linear}. Finally, in Section~\ref{IP nonlinear} we prove with the help of  a suitable integral identity (Lemma~\ref{lem:integral identity}) and the linearizatzion of the DN map the main theorem on the inverse problem for the nonlocal viscous wave equation with a nonlinear perturbation (Theorem~\ref{Thm: main nonlinear}).
	
	\section{Preliminaries}
	\label{sec: preliminaries}
	
    In this section, we introduce several function spaces together with the fractional Laplacian, recall some important properties of this nonlocal operator and describe an abstract framework for solving some classes of second order in time PDEs.
	
	\subsection{Fractional Sobolev spaces and fractional Laplacian}
	\label{sec: Fractional Sobolev spaces}
	
	We denote by $\schwartz(\R^n)$ and $\tempered(\R^n)$ Schwartz functions and tempered distributions respectively. We define the Fourier transform by
	\begin{equation}
		\fourier u(\xi) \vcentcolon = \hat{u}(\xi)\vcentcolon = \int_{\R^n} u(x)e^{-\mathrm{i}x \cdot \xi} \,dx.
	\end{equation}
	By duality it can be extended to the space of tempered distributions and will again be denoted by $\fourier u = \widehat{u}$, where $u \in \tempered(\R^n)$, and we denote the inverse Fourier transform by $\ifourier$. 
	
	Given $s\in\R$, the fractional Sobolev space $H^{s}(\R^{n})$ is the set of all tempered distributions $u\in\tempered(\R^n)$ such that
	\begin{equation}\notag
		\|u\|_{H^{s}(\mathbb{R}^{n})}\vcentcolon = \left\|\langle D\rangle^s u \right\|_{L^2(\R^n)}<\infty,
	\end{equation}
	where $\langle D\rangle^s$ is the Bessel potential operator of order $s$ having Fourier symbol $\LC 1+|\xi|^2\RC^{s/2}$. 
	The fractional Laplacian of order $s\geq 0$ can be defined as a Fourier multiplier
	\begin{equation}\label{eq:fracLapFourDef}
		(-\Delta)^{s} u = \ifourier
  (\abs{\xi}^{2s}\widehat{u}(\xi)),
	\end{equation}
	for $u \in \tempered(\R^n)$ whenever the right hand side of the above identity is well-defined. 
	In addition, it is also known that for $s\geq 0$, an equivalent norm on $H^s(\R^n)$ is given by
	\begin{equation}
		\label{eq: equivalent norm on Hs}
		\|u\|_{H^s(\R^n)}^*= \|u\|_{L^{2}(\mathbb{R}^{n})}+\|(-\Delta)^{s/2}u  \|_{L^{2}(\mathbb{R}^{n})},
	\end{equation}
	and the fractional Laplacian $(-\Delta)^{s}\colon H^{t}(\R^n) \to H^{t-2s}(\R^n)$ is a bounded linear operator for all $s \geq0$ and $t \in \R$. 
    
    The above heuristically introduced UCP reads more formally as follows: 
	
	\begin{proposition}[{UCP for fractional Laplacians}]\label{prop:UCP} 
		Let $s> 0$ be a non-integer and $t\in\R$. If $u\in H^t(\R^n)$ satisfies $u=(-\Delta)^s u=0$ in a nonempty open subset $V\subset \R^n$, then $u\equiv 0$ in $\R^n$.
	\end{proposition}
	
	The preceding proposition was first shown in \cite[Theorem 1.2]{GSU20} for the range $s\in (0,1)$, in which case the fractional Laplacian $(-\Delta)^s$ can be equivalently computed as the singular integral 
	\[
	(-\Delta)^su(x)=C_{n,s}\text{p.v.}\int_{\R^n}\frac{u(x)-u(y)}{|x-y|^{n+2s}}\,dy
	\]
	for sufficiently nice functions $u$ and some constant $C_{n,s}>0$. For the higher order case $s>1$, one can apply the standard Laplacian to the equation, then the classical UCP for the Laplacian yields iteratively the desired result.
	
 For the well-posedness theory, we will use the following Poincar\'e inequality.
	\begin{proposition}[{Poincar\'e inequality (cf.~\cite[Lemma~5.4]{RZ-unbounded})}]\label{prop:Poincare ineq} Let $\Omega\subset\R^n$ be a bounded domain. For any $s\geq 0$, there exists $C>0$ such that
		\begin{equation}
			\label{eq: poincare ineq}
			\|u\|_{L^2(\Omega)}\leq C  \|(-\Delta)^{s/2}u  \|_{L^2(\R^n)}
		\end{equation}
		for all $u\in C_c^{\infty}(\Omega)$.
	\end{proposition}
	
	Next we introduce some local variants of the above fractional Sobolev spaces. If $\Omega\subset \R^n$ is an open set, $F\subset \R^n$ a closed set and $s\in\mathbb{R}$,
	then we set
	\begin{align*}
		H^{s}(\Omega) & \vcentcolon =\left\{ u|_{\Omega}\,; \, u\in H^{s}(\mathbb{R}^{n})\right\},\\
		\widetilde{H}^{s}(\Omega) & \vcentcolon =\text{closure of \ensuremath{C_{c}^{\infty}(\Omega)} in \ensuremath{H^{s}(\mathbb{R}^{n})}},\\
		H_F^s&\vcentcolon =\{u\in H^s(\R^n)\,;\,\supp(u)\subset F\}.
	\end{align*}
	Meanwhile, $H^{s}(\Omega)$ is a Banach space with respect to the quotient norm
	\[
	\|u\|_{H^{s}(\Omega)}\vcentcolon =\inf\left\{ \|U\|_{H^{s}(\mathbb{R}^{n})} \,; \,  U\in H^{s}(\mathbb{R}^{n})\mbox{ and }U|_{\Omega}=u\right\} .
	\]
	Hence, using the fact that \eqref{eq: equivalent norm on Hs} is an equivalent norm on $\widetilde{H}^s(\Omega)$, Propositions~\ref{prop:Poincare ineq} and the density of $C_c^{\infty}(\Omega)$ in $\widetilde{H}^s(\Omega)$, we have:
	\begin{lemma}
		\label{lemma: equivalent norm on tilde spaces}
		Let $\Omega\subset\R^n$ be a bounded domain and $s\geq 0$. Then an equivalent norm on $\widetilde{H}^s(\Omega)$ is given by
		\begin{equation}
			\label{eq: equivalent norm on tilde spaces}
			\|u\|_{\widetilde{H}^s(\Omega)}=\|(-\Delta)^{s/2}u\|_{L^2(\R^n)}.
		\end{equation}
	\end{lemma}
	The observation of Lemma~\ref{lemma: equivalent norm on tilde spaces} will be of constant use in the well-posedness theory below.
	
	\subsection{Bochner spaces}
	\label{subsec: Bochner spaces}
	
	Next, we introduce some standard function spaces for time-dependent PDEs adapted to the nonlocal setting considered in this article. Let $X$ be a Banach space and $(a,b)\subset\R$. Then we let $C^k([a,b]\,;X)$, $L^p(a,b\,;X)$ ($k\in\N,1\leq p\leq \infty$) stand for the space of $k-$times continuously differentiable functions and the space of measurable functions $u\colon (a,b)\to X$ such that $t\mapsto \|u(t)\|_X\in L^p([a,b])$. These spaces carry the norms
	\begin{equation}
		\label{eq: Bochner spaces}
		\begin{split}
			\|u\|_{L^p(a,b\,;X)}&\vcentcolon = \left(\int_a^b\|u(t)\|_{X}^p\,dt\right)^{1/p}<\infty,\\
			\|u\|_{C^k([a,b];X)}&\vcentcolon = \sup_{0\leq \ell\leq k}\|\partial_t^{\ell}u\|_{L^{\infty}([a,b];X)}
		\end{split}
	\end{equation} 
    with the usual modifications in the case $p=\infty$.
	
	Additionally, whenever $u\in L^1_{\loc}(a,b\,;X)$ with $X$ being a space of functions over a subset of some euclidean space, such as $L^2(\Omega)$ or $H^s(\R^n)$, then $u$ is identified with a function $u(x,t)$ and $u(t)$ denotes the function $ x\mapsto u(x,t)$ for almost all $t$. This is justified by the fact, that any $u\in L^q(a,b\,;L^p(\Omega))$ with $1\leq q,p<\infty$ can be seen as a measurable function $u\colon \Omega\times (a,b)\to \R$ such that the norm $\|u\|_{L^q(a,b\,;L^p(\Omega))}$, as defined in \eqref{eq: Bochner spaces}, is finite. In particular, one has $L^p(0,T;L^p(\Omega))=L^p(\Omega_T)$ for $1\leq p<\infty$. Clearly, a similar statement holds for the spaces $L^q(a,b\,;H^{s}(\R^n))$ and their local versions. If no confusion arises we also denote $L^p(0,T;X)$ by $L^p(X)$ and $L^q(0,T;L^p(\Omega))$ by $L^qL^p$.
 
    Furthermore, the distributional derivative $\frac{du}{dt}\in \distr((a,b)\,;X)$ is identified with the derivative $\partial_tu\in \distr(\Omega\times (a,b))$ as long as it is well-defined. Here $ \distr((a,b)\,;X)$ stands for all continuous linear operators from $C_c^{\infty}((a,b))$ to $X$.

    \subsection{Abstract framework for solving some 2nd order in time PDEs}
    \label{subsec: abstract framework}

    We collect here preliminary material to solve  some classes of second order PDEs and for a more comprehensive presentation the interested reader can consult \cite[Chapter XVIII]{DautrayLionsVol5}. 

    \begin{definition}
    \label{def: usual conditions}
        We say that a tuple $(V,H,a,b,c)$ consisting of two complex Hilbert spaces $V,H$, two families of sesquilinear forms $a,b$ over $V$ and one family of sesquilinear forms $c$ over $H$ satisfy the \emph{usual conditions} if they have the following properties: The spaces $V$ and $H$ are such that 
    \begin{equation}
    \label{eq: Hilbert triplet condition}
        V\hookrightarrow H\hookrightarrow V'
    \end{equation}
    and the inclusions are dense.\footnote{Here, $V'$ denotes the antidual of $V$, that is the space of all antilinear continuous functionals on $V$, and the antidual of $H$ is identified with $H$ via Riesz's representation theorem. Moreover, the dual of $V$, $H$ are still denoted by $V^*$ and $H^*$. Recall that for any Hilbert space $X$, we have $X'=\overline{X^*}$.} The triple $(a,b,c)$ fulfill the following assumptions:
    \begin{enumerate}
        \item[(S1)]\label{S1} $a(t;\cdot,\cdot)$, $t\in [0,T]$, is a family of sesquilinear forms over $V$, which can be decomposed as 
        \begin{equation}
        \label{eq: decomp of a}
            a=a_0+a_1.
        \end{equation}
        The principal part $a_0$ is required to satisfy
        \begin{enumerate}
            \item[(A1)]\label{A1} $a_0(t;u,v)\in C^1([0,T])$ for all $u,v\in V$ such that the sesquilinear forms $a_0,\partial_t a_0$ are continuous over $V$,
            \item[(A2)]\label{A2} $a_0$ is Hermitian (i.e. $a_0(t;u,v)=\overline{a_0(t;v,u)}$),
            \item[(A3)]\label{A3} $a_0$ is coercive over $V$ with respect to $H$ in the sense that
            \begin{equation}
            \label{eq: weakly coercive}
                a_0(t;u,u)\geq \alpha\|u\|_V^2-\lambda \|u\|_H^2
            \end{equation}
            for some $\alpha>0,\lambda\in\R$ and all $u\in V$, $t\in [0,T]$.
        \end{enumerate}
        The lower order part $a_1$ is assumed to satisfy
        \begin{enumerate}
            \item[(A4)]\label{A4} $a_1(t;u,v)\in C([0,T])$ for all $u,v\in V$,
            \item[(A5)]\label{A5} $|a_1(t;u,v)|\leq C\|u\|_V\|v\|_H$ for all $u,v\in V$ and $t\in [0,T]$.
        \end{enumerate}
        \item[(S2)]\label{S2} $b(t;\cdot,\cdot)$, $t\in [0,T]$, is a family of sesquilinear forms over $V$, which can be decomposed as 
        \begin{equation}
        \label{eq: decomp of b}
            b=b_0+b_1.
        \end{equation}
        \begin{enumerate}
            \item[(B1)]\label{B1} $b_0$ is a continuous sesquilinear form over $V$, Hermitian and coercive in the sense 
            \begin{equation}
                b_0(t;u,u)\geq \mu\|u\|_V^2
            \end{equation}
            for some $\mu>0$ and $t\in [0,T],v\in V$,
            \item[(B2)]\label{B2} $b_1$ is a sesquilinear form satisfying
            \begin{equation}
            \label{eq: continuity of b1}
                |b_1(t;u,v)|\leq C\|u\|_V\|v\|_H
            \end{equation}
            for all $t\in [0,T]$ and $u,v\in V$,
            \item[(B3)]\label{B3} $b_j(t;u,v)\in C([0,T])$ for all $u,v\in V$ and $j=0,1$.
        \end{enumerate}
        \item[(S3)]\label{S3} $c(t;\cdot,\cdot)$, $t\in [0,T]$, is a family of sesquilinear forms over $H$, which can be written as
        \begin{equation}
        \label{eq: def of c}
            c(t;u,v)=\langle C(t)u,v\rangle_H
        \end{equation}
        for $t\in [0,T]$ and $u,v\in H$. Here $C(t)$, $t\in [0,T]$, is a family of linear bounded operators on $H$ to itself such that
        \begin{enumerate}
            \item[(C1)]\label{C1} $C(t)$ is Hermitian and coercive over $H$, that is
            \begin{equation}
                \langle C(t)u,u\rangle_H\geq \gamma \|u\|^2_H
            \end{equation}
            for some $\gamma>0$ and all $u\in H$, $t\in[0,T]$,
            \item[(C2)]\label{C2} $\langle C(t)u,v\rangle_H\in C^1([0,T])$ for all $u,v\in H$.
        \end{enumerate}
    \end{enumerate}
    \end{definition}

     \begin{example}
    \label{example spaces}
         Let $\Omega\subset\R^n$ be any bounded Lipschitz domain. Then the separable Hilbert spaces $V=\widetilde{H}^s(\Omega)$ and $H=L^2(\Omega)$\footnote{These spaces are considered here as being complex, but later on we always assume that they consist of real valued functions.} satisfy \eqref{eq: Hilbert triplet condition} and the inclusions are dense. This is a direct consequence of the Sobolev embedding, the assumption that $\Omega$ is a bounded Lipschitz domain and the fact that $u\in\widetilde{H}^s(\Omega)$ implies $u=0$ a.e. in $\Omega^c$. 
    \end{example}
   
    \begin{example}
    \label{example sesquilinear forms}
        Let $V,H$ be as in Example~\ref{example spaces}. Then we define the principle part sesquilinear form $a_0\colon V\times V\to\C$
        \begin{equation}
        \label{eq: def a0}
            a_0(u,v)=\langle (-\Delta)^{s/2} u,(-\Delta)^{s/2} v\rangle_{L^2(\R^n)}
        \end{equation}
       for all $u,v\in V$. This form clearly satisfies (A1)-(A2) (see Lemma~\ref{lemma: equivalent norm on tilde spaces}). The condition (A3) is a consequence of the Poincar\'e inequality (Proposition~\ref{prop:Poincare ineq}). More precisely, there holds
       \begin{equation}
       \label{eq: coercivity example a0}
           a_0(u,u)\geq c\|u\|_V^2.
       \end{equation}
       Hence, by the properties of $a_0$, particularly \eqref{eq: coercivity example a0}, we can choose $b_0=a_0$ and the properties (B1) and (B3) are fulfilled. Moreover, we set $b_1=0$ and $C(t)=\text{id}_{L^2(\Omega)}$. Finally, we assume that $q\in L^1_{loc}(\Omega_T)$ is any function such that the induced sesquilinear forms $a_1(t;\cdot,\cdot)\colon V\times V\to\C$, $t\in [0,T]$, given by
       \begin{equation}
       \label{eq: potential term}
           a_1(t;u,v)=\langle q(t)u,v\rangle_{L^2(\Omega)}
       \end{equation}
       are well-defined and satisfy (A4), (A5). Hence, the tuple $(V,H,a,b,c)$ satisfies the usual conditions in the sense of Defintion~\ref{def: usual conditions}.
\end{example}

    Next, we introduce several function spaces used throughout this article. 
    \begin{definition}
         Let $T>0$. Suppose that we have given Hilbert spaces $V,H$ satisfying the conditions in Definition~\ref{def: usual conditions} and a family of bounded linear operators $C(t)\in L(H)$, $t\in[0,T]$, such that the related sesquilinear forms
         \[
            c(t;u,v)=\langle C(t)u,v\rangle_H
         \]
         fulfill the property (S3) of Definition~\ref{def: usual conditions}. 
    \begin{enumerate}
        \item[(F1)] Then we set
        \begin{equation}
        \label{eq: space W(V)}
            W_c(0,T;V)=\{v\in L^2(0,T;V)\,;\,\frac{d}{dt}(Cv)\in L^2(0,T;V')\},
        \end{equation}
        which carries the norm
        \begin{equation}
        \label{eq: norm W(V)}
            \|v\|_{W_c(0,T;V)}=\left(\|v\|_{L^2(0,T;V)}^2+\|\frac{d}{dt}Cv\|_{L^2(0,T;V')}^2\right)^{1/2}.
        \end{equation}
        \item[(F2)] Furthermore, we define
        \begin{equation}
        \label{eq: solution space}
            \widetilde{W}_c(0,T;V)=\{u\in L^2(0,T;V)\,;\,\partial_t u\in W_c(0,T;V)\}
        \end{equation}
        and equip it with the norm
        \begin{equation}
        \label{eq: norm solution space}
            \|u\|_{\widetilde{W}_c(0,T;V)}=\left(\|u\|_{L^2(0,T;V)}^2+\|\partial_t u\|_{W_c(0,T;V)}^2\right)^{1/2}
        \end{equation}
    \end{enumerate}
    \end{definition}
    \begin{remark}
        Clearly both space $W_c(0,T;V)$ and $\widetilde{W}_c(0,T;V)$ are Hilbert spaces. Moreover, if $C(t)=\text{id}_H$, then we drop the subscript $c$.
    \end{remark}
    The next lemma collects a few properties of these spaces (see \cite[Chapter XVIII, \S 5]{DautrayLionsVol5}).
    \begin{lemma}
    \label{lemma: embeddings}
         Let $T>0$. Suppose that we have given Hilbert spaces $V,H$ satisfying the conditions in Definition~\ref{def: usual conditions} and a family of bounded linear operators $C(t)\in L(H)$, $t\in[0,T]$, such that the related sesquilinear forms
         \[
            c(t;u,v)=\langle C(t)u,v\rangle_H
         \]
         fulfill the property (S3) of Definition~\ref{def: usual conditions}. 
        \begin{enumerate}[(i)]
            \item One has the embeddings
        \begin{equation}
        \label{eq: continuous embeddings}
            \widetilde{W}_c(0,T;V)\hookrightarrow C([0,T];V)\text{ and }W_c(0,T;V)\hookrightarrow C([0,T];H).
        \end{equation}
        \item The space $C_c^{\infty}([0,T];V)$ is dense in $W_c(0,T;V)$ and in $\widetilde{W}_c(0,T;V)$.
        \end{enumerate}
    \end{lemma}
    
    Now, we can formulate the abstract forward problem.
    \begin{problem}
    \label{linear problem}
        Suppose we have given a tuple $(V,H,a,b,c)$ satisfying the usual conditions. Does there exist for all functions $u_0\in V$, $u_1\in H$ and $f\in L^2(0,T;V')$ a unique function $u\in \widetilde{W}_c(0,T;V)$ satisfying
        \begin{equation}
        \label{eq: general parabolic PDE}
            \frac{d}{dt}\,c(\cdot;\partial_t u,v)+b(\cdot;\partial_t u,v)+a(\cdot;u,v)=\langle f,v\rangle_{V'\times V}
        \end{equation}
        for all $v\in V$ in the sense of $\distr((0,T))$
        and
        \begin{equation}
        \label{eq: initial conditions}
            u(0)=u_0\text{ in }V\quad\text{ and }\quad \partial_t u(0)=u_1\text{ in }H.
        \end{equation}
    \end{problem}
    \begin{remark}
        Note that $u\in \widetilde{W}_c(0,T;V)$ guarantees that all terms in  \eqref{eq: general parabolic PDE} are well-defined and the embeddings \eqref{eq: continuous embeddings} ensure that  \eqref{eq: initial conditions} makes sense.
    \end{remark}

    Now, we can state the main well-posedness result on abstract second order in time PDEs, which we will use later on.

    \begin{theorem}[Well-posedness abstract PDEs]
    \label{thm: abstract well-posedness}
        Let $T>0$. Assume $(V,H,a,b,c)$ consisting of two complex Hilbert spaces $V,H$, two sesquilinear forms $a,b$ over $V$ and a sesquilinear form over $H$ satisfy the usual conditions. Then for any $u_0\in V$, $u_1\in H$ and $f\in L^2(0,T;V')$, the Problem~\ref{linear problem} has a unique solution $u\in \widetilde{W}_c(0,T;V)$.
    \end{theorem}
    \begin{proof}
        This result is a direct consequence of \cite[Chapter XVIII, \S 5, Theorem~1, Remark~4]{DautrayLionsVol5}.
    \end{proof}
	
	\section{Well-posedness theory of viscous wave equations}
	\label{sec: well-posedness}
     In this section, we study the well-posedness theory of the viscous wave equation with linear and nonlinear perturbations.
	
	\subsection{Viscous wave equation with linear perturbations}
    \label{subsec: well-posedness linear}
	
	Let us start by stating the well-posedness result in the linear case.
	
	\begin{theorem}
		\label{thm: well-posedness linear}
		Let $\Omega\subset\R^n$ be a bounded Lipschitz domain, $T>0$ and $s>0$. Suppose that the (real valued) function $q\in L^1_{loc}(\Omega_T)$ has the following properties:
        \begin{enumerate}[(i)]
        \item\label{integrability of q} $q\in L^{\infty}(0,T;L^p(\Omega))$ for some $1\leq p<\infty$ satisfying
        \begin{equation}
        \label{eq: restriction of p in well-posedness linear}
		\begin{cases}
			n/s\leq p\leq \infty, &\, \text{if }\, 2s< n,\\
			2<p\leq \infty,  &\, \text{if }\, 2s= n,\\
			2\leq p\leq \infty, &\, \text{if }\, 2s > n,
		\end{cases}
		\end{equation}
            \item\label{continuity of q} $t\mapsto \int_{\Omega}q(x,t)\varphi(x)\,dx\in C([0,T])$ for any $\varphi\in C_c^{\infty}(\Omega)$.
        \end{enumerate}
		Then for any pair $(u_0, u_1)\in \widetilde{H}^s(\Omega)\times L^2(\Omega)$ and $h\in L^2(0,T;H^{-s}(\Omega))$ there exists a unique solution $u\in \widetilde{W}(0,T;\widetilde{H}^s(\Omega))$ of
		\begin{equation}
			\label{eq: well-posedness linear case}
			\begin{cases}
				\LC \partial_t^2 +(-\Delta)^s\partial_t +(-\Delta)^s+q\RC u=h &\text{ in }\Omega_T\\
				u=0     &\text{ in }(\Omega_e)_T,\\
				u(0)=u_0,\quad \partial_t u(0)=u_1  &\text{ in }\Omega.
			\end{cases}
		\end{equation}
        Moreover, $u$ satisfies the following energy identity\footnote{Here, and throughout this work, we use the notation $\langle \cdot , \cdot \rangle$ for the natural duality pairing between the corresponding function spaces. In the present context, it denotes the pairing between $H^{-s}(\Omega)$ and $\widetilde{H}^s(\Omega)$, while later it is also used for the pairing between $L^2(0,T;H^{-s}(\Omega))$ and $L^2(0,T;\widetilde{H}^s(\Omega))$.}
		\begin{equation}
			\label{eq: energy identity linear case}
			\begin{split}
				& \|\partial_t u(t) \|_{L^2(\Omega)}^2+ \|(-\Delta)^{s/2}u(t)\|_{L^2(\R^n)}^2+ 2\|(-\Delta)^{s/2}\partial_tu\|_{L^2(\R^n_t)}^2\\
				&=\|u_1 \|_{L^2(\Omega)}^2+\|(-\Delta)^{s/2}u_0 \|_{L^2(\R^n)}^2+2\int_0^t\langle h(\tau),\partial_t u (\tau)\rangle\,d\tau-2\langle q u,\partial_t u \rangle_{L^2(\Omega_t)}
			\end{split}
		\end{equation}
		for all $t\in [0,T]$. Moreover, if $(u_{0,j},u_{1,j})\in \widetilde{H}^s(\Omega)\times L^2(\Omega)$, $h_j\in L^2(0,T;H^{-s}(\Omega))$ and $u_j\in \widetilde{W}(0,T;\widetilde{H}^s(\Omega))$ denote the related unique solution to \eqref{eq: well-posedness linear case} for $j=1,2$, then the following continuity estimate holds
		\begin{equation}
			\label{eq: continuity estimate}
			\begin{split}
				&\|u_1-u_2\|_{L^{\infty}(0,T;\widetilde{H}^s(\Omega))}+\|\partial_t u_1-\partial_t u_2\|_{L^{\infty}(0,T;L^2(\Omega))}+\|\partial_t u_1-\partial_t u_2\|_{L^{2}(0,T;\widetilde{H}^s(\Omega))}\\
				&\quad \leq C(\|u_{0,1}-u_{0,2}\|_{\widetilde{H}^s(\Omega)}+\|u_{1,1}-u_{1,2}\|_{L^2(\Omega)}+\|h_1-h_2\|_{L^2(0,T;H^{-s}(\Omega))})
			\end{split}
		\end{equation}
		for some $C>0$ depending on $T>0$.
	\end{theorem}
	
	\begin{proof}
		For the time being assume $\widetilde{H}^s(\Omega)$ and $L^2(\Omega)$ consist of complex functions. We claim that we are in the setting of Example~\ref{example sesquilinear forms}. To this end, we only need to verify that the sesquilinear form 
        \begin{equation}
        \label{eq: sesquilinear form from potential}
            a_1(t;u,v)=\langle q(t)u,v\rangle_{L^2(\Omega)}
        \end{equation}
        for $u,v\in\widetilde{H}^s(\Omega)$ satisfies (A4) and (A5). If we can show the estimate
		\begin{equation}
			\label{eq: estimate a1}
        \begin{split}
            \left|\langle qu,v\rangle_{L^2(\Omega)} \right|&\leq C\|q(t)\|_{L^p(\Omega)}\|u\|_{\widetilde{H}^s(\Omega)}\|v\|_{L^2(\Omega)}\\
            &\leq C\|q\|_{L^{\infty}(0,T;L^p(\Omega))}\|u\|_{\widetilde{H}^s(\Omega)}\|v\|_{L^2(\Omega)}
        \end{split}
		\end{equation}
		for some $C>0$ independent of $t\in [0,T]$, then (A5) follows. The case $p=\infty$ is clear.  In the case $\frac{n}{s}\leq p<\infty$ with $2s< n$ one can use H\"older's inequality with
		\[
		\frac{1}{2}=\frac{n-2s}{2n}+\frac{s}{n},
		\]
		$L^{r_2}(\Omega)\hookrightarrow L^{r_1}(\Omega)$ for $r_1\leq r_2$ as $\Omega\subset\R^n$ is bounded and Sobolev's inequality to obtain
		\begin{equation}
			\label{eq: computation for L2 estimate}
			\begin{split}
				\left|\langle qu,v\rangle_{L^2(\Omega)} \right|&\leq \|qu\|_{L^2(\Omega)}\|v\|_{L^2(\Omega)}\\
				&\leq \|q\|_{L^{n/s}(\Omega)}\|u\|_{L^{\frac{2n}{n-2s}}(\Omega)}\|v\|_{L^2(\Omega)}\\
				&\leq C\|q\|_{L^{n/s}(\Omega)}\|u\|_{L^{\frac{2n}{n-2s}}(\Omega)}\|v\|_{L^2(\Omega)}\\
				&\leq C\|q\|_{L^{p}(\Omega)}\|(-\Delta)^{s/2} u \|_{L^2(\R^n)}\|v\|_{L^2(\Omega)}
			\end{split}
		\end{equation}
		for all $u,v\in \widetilde{H}^s(\Omega)$. In the case $2s>n$ one can use the embedding $H^s(\R^n)\hookrightarrow L^{\infty}(\R^n)$ together with Lemma~\ref{lemma: equivalent norm on tilde spaces} and the boundedness of $\Omega$ to see that the estimate \eqref{eq: computation for L2 estimate} holds. In the case $n=2s$ one can use the boundedness of the embedding $\widetilde{H}^s(\Omega)\hookrightarrow L^{\overline{p}}(\Omega)$ for all $2\leq \overline{p}<\infty$, H\"older's inequality and the boundedness of $\Omega$ to get the estimate \eqref{eq: computation for L2 estimate}. In fact, the aforementioned embedding in the critical case follows by \cite{Ozawa} and the Poincar\'e inequality. 

        Next, we show that the condition \ref{continuity of q} implies (A4). For this purpose, let us define for any $t\in [0,T]$ the function
        \[
            \Phi_{\varphi}(t)=\int_{\Omega}q(x,t)\varphi(x)\,dx
        \]
        for any $\varphi\in L^1_{loc}(\Omega)$ such that the integral is well-defined. Suppose $u_j,v_j\in \widetilde{H}^s(\Omega)$, $j=1,2$, are given. Then the estimate \eqref{eq: computation for L2 estimate} implies
        \begin{equation}
        \label{eq: continuity estimate}
        \begin{split}
            &|\Phi_{u_1\overline{v_1}}(t)-\Phi_{u_2\overline{v_2}}(t)|\\
            &\leq |\Phi_{(u_1-u_2)\overline{(v_1-v_2)}}(t)|+|\Phi_{(u_1-u_2)\overline{v_2}}(t)|+|\Phi_{u_2\overline{(v_1-v_2)}}(t)|\\
            &\leq C(\|u_1-u_2\|_{\widetilde{H}^s(\Omega)}\|v_1-v_2\|_{L^2(\Omega)}+\|u_1-u_2\|_{\widetilde{H}^s(\Omega)}\|v_2\|_{L^2(\Omega)}+\|u_2\|_{\widetilde{H}^s(\Omega)})
            \end{split}
        \end{equation}
        for all $t\in [0,T]$. Now, for fixed $u,v\in \widetilde{H}^s(\Omega)$ there exist sequences $u_k,v_k\in C_c^{\infty}(\Omega)$ such that 
        \[
        u_k\to u\text{ and }v_k\to v\text{ in }\widetilde{H}^s(\Omega)
        \]
        as $k\to\infty$. Then the estimate \eqref{eq: continuity estimate} shows that
        \[
            \Phi_{\varphi_k}(t)\to \Phi_{u\overline{v}}(t)\text{ as }k\to\infty
        \]
       uniformly in $t\in[0,T]$, where $\varphi_k=u_k\overline{v_k}\in C_c^{\infty}(\Omega)$. Hence, for any $\eps>0$ there exists $k_0\in\N$ such that 
       \[
            \|\Phi_{\varphi_k}-\Phi_{u\overline{v}}\|_{L^{\infty}([0,T])}<\eps/3
       \]
       for all $k\geq k_0$. Now, let us fix such a $k\geq k_0$. On the other hand, for any $k\in\N$ we know by \ref{continuity of q} that $\Phi_{\varphi_k}\in C([0,T])$ and thus for given $k\geq k_0$ we can choose $\delta=\delta(k)>0$ such that if $t_1,t_2\in [0,T]$ with $|t_1-t_2|<\delta$, we have
       \[
            |\Phi_{\varphi_k}(t_1)-\Phi_{\varphi_k}(t_2)|<\eps/3.
       \]
       Thus, we get 
       \[
       \begin{split}
           &|\Phi_{u\overline{v}}(t_1)-\Phi_{u\overline{v}}(t_2)|\\
           &\leq |\Phi_{u\overline{v}}(t_1)-\Phi_{\varphi_k}(t_1)|+|\Phi_{u\overline{v}}(t_2)-\Phi_{\varphi_k}(t_2)|+|\Phi_{\varphi_k}(t_1)-\Phi_{\varphi_k}(t_2)|<\eps
       \end{split}
       \]
       for all $t_1,t_2\in [0,T]$ such that $|t_1-t_2|<\delta$. Hence, we have $\Phi_{u\overline{v}}\in C([0,T])$ and so the sesquilinear form $a_1$ satisfies (A4) as well. 

       Therefore, if we extend $h\in L^2(0,T;H^{-s}(\Omega))=L^2(0,T;(\widetilde{H}^s(\Omega))^*)$ to its unique antilinear functional $H\in L^2(0,T;(\widetilde{H}^s(\Omega))')$, we may deduce from Theorem~\ref{thm: abstract well-posedness} the existence of a unique solution $\widetilde{u}\in\widetilde{W}(0,T;\widetilde{H}^s(\Omega))$ of \eqref{eq: general parabolic PDE}, where $(V,H,a,b,c)$ are as above. Since $q$ is real valued, one easily deduces the existence of a unique real valued solution $u\in \widetilde{W}(0,T;\widetilde{H}^s(\Omega))$ of \eqref{eq: well-posedness linear case}.

       The energy identity \eqref{eq: energy identity linear case} and the continuity estimate \eqref{eq: continuity estimate} are direct consequences of \cite[Chapter XVIII, \S 5, Equations (5.81)-(5.82) and Remark 4]{DautrayLionsVol5} and the fact that all appearing functions are real valued.
	\end{proof}

    Next, we introduce several function spaces, which are used below. First, we define 
	\begin{equation}
		\label{eq: nice solutions}
		\widetilde{W}_*(0,T;H^s(\R^n))=\{v\in \widetilde{W}(0,T;H^s(\R^n))\,;\,v(0)\in \widetilde{H}^s(\Omega),\,\partial_t v(0)\in L^2(\Omega)\}.
	\end{equation}
	By Lemma~\ref{lemma: embeddings} it follows that $\widetilde{W}_*(0,T;H^s(\R^n))\subset \widetilde{W}(0,T;H^s(\R^n))$ is a closed subspace and thus again a Hilbert space. Additionally, we introduce the Banach space
    \begin{equation}
    \label{eq: solution space non zero exterior conditions}
        \widetilde{W}_{ext}(0,T;\widetilde{H}^s(\Omega))=\widetilde{W}(0,T;\widetilde{H}^s(\Omega))+\widetilde{W}_*(0,T;H^s(\R^n))\subset H^1(0,T;H^s(\R^n)),
    \end{equation}
    which is as usual endowed with the norm
   \begin{equation}
   \label{eq: norm on sum}
       \|u\|_{\widetilde{W}_{ext}(0,T;\widetilde{H}^s(\Omega))}=\inf(\|v\|_{\widetilde{W}(0,T;\widetilde{H}^s(\Omega))}+\|\varphi\|_{\widetilde{W}(0,T;H^s(\R^n))}),
   \end{equation}
   where the infimum is taken over all $v\in \widetilde{W}(0,T;\widetilde{H}^s(\Omega))$ and $\varphi\in \widetilde{W}_*(0,T;H^s(\R^n))$ such that $u=v+\varphi$. Later $\widetilde{W}_{ext}(0,T;\widetilde{H}^s(\Omega))$ will play the role of the solution space to problems with nonzero exterior conditions $\varphi$.
    
    The last space we need is
	\[
	\widetilde{W}^s_\text{rest}((\Omega_e)_T)=\{v|_{(\Omega_e)_T}\,;\,v\in \widetilde{W}_*(0,T;H^s(\R^n))\},
	\]
	which is a Banach space when endowed with the corresponding quotient norm.
    
    \begin{corollary}
    \label{Well-posedness with nonzero exterior conditions}
        Let $\Omega\subset\R^n$ be a bounded Lipschitz domain, $T>0$ and $s>0$. Suppose that the (real valued) function $q\in L^1_{loc}(\Omega_T)$ satisfies the conditions in Theorem~\ref{thm: well-posedness linear}.
        \begin{enumerate}[(i)]
        \item\label{item: uniqueness linear} For any pair $(u_0, u_1)\in \widetilde{H}^s(\Omega)\times L^2(\Omega)$, $h\in L^2(0,T;H^{-s}(\Omega))$ and $\Phi\in \widetilde{W}_*(0,T;\widetilde{H}^s(\R^n))$, there exists a unique solution $u\in \widetilde{W}_{ext}(0,T;\widetilde{H}^s(\Omega))$ of
		\begin{equation}
			\label{eq: well-posedness linear case nonzero exterior cond}
			\begin{cases}
				\LC \partial_t^2 +(-\Delta)^s\partial_t +(-\Delta)^s+q\RC u=h &\text{ in }\Omega_T\\
				u=\Phi     &\text{ in }(\Omega_e)_T,\\
				u(0)=u_0,\quad \partial_t u(0)=u_1  &\text{ in }\Omega.
			\end{cases}
		\end{equation}
        This means
        \begin{enumerate}[(I)]
            \item\label{linear PDE with ext cond} for all $v\in\widetilde{H}^s(\Omega)$ one has
            \begin{equation}
            \label{eq: weak solution linear nonzero ext cond}
            \begin{split}
                 &\frac{d}{dt}\,\langle \partial_t u,v\rangle_{L^2(\Omega)}+\langle (-\Delta)^{s/2}\partial_t u,(-\Delta)^{s/2}v\rangle_{L^2(\R^n)}\\
                 &+\langle (-\Delta)^{s/2}u,(-\Delta)^{s/2}v\rangle_{L^2(\R^n)}+\langle q u,v\rangle_{L^2(\Omega)}=\langle h,v\rangle
            \end{split}
        \end{equation}
         in the sense of $\distr((0,T))$,
            \item\label{exterior cond} $u=\Phi$ in $(\Omega_e)_T$,
            \item\label{initial values} $u(0)=u_0$ in $\widetilde{H}^s(\Omega)$ and $\partial_t u(0)=u_1$ in $L^2(\Omega)$.
        \end{enumerate}
        \item\label{item: sol indep of extension} If $\varphi\in \widetilde{W}^s_{rest}((\Omega_e)_T)$ and $\Phi_1,\Phi_2$ are any two representations of $\varphi$ with unique solutions $u_1$ and $u_2$ of \eqref{eq: well-posedness linear case nonzero exterior cond} with $\Phi=\Phi_1$ and $\Phi=\Phi_2$, respectively, then there holds $u_1=u_2$. In particular, for any $\varphi\in \widetilde{W}^s_{rest}((\Omega_e)_T)$ we have a unique solution $u$ of \eqref{eq: well-posedness linear case nonzero exterior cond}.
        \end{enumerate}
    \end{corollary}
    \begin{proof}
        \ref{item: uniqueness linear}: We first observe that if $u\in\widetilde{W}_{ext}(0,T;\widetilde{H}^s(\Omega))$ and $\Phi\in \widetilde{W}_*(0,T;H^s(\R^n))$, then one has $u=\Phi$ in $(\Omega_e)_T$ if and only if  
        \begin{equation}
        \label{eq: equivalent characterization of ext cond}
            u-\Phi\in \widetilde{W}(0,T;\widetilde{H}^s(\Omega)).
        \end{equation}
        In fact, if $u=v+\psi$ with $v\in \widetilde{W}(0,T;\widetilde{H}^s(\Omega))$ and $\psi\in \widetilde{W}_*(0,T;H^s(\R^n))$, then one has $u=\Phi$ in $(\Omega_e)_T$ if and only if $\psi=\Phi$ in $(\Omega_e)_T$. As $\Omega$ has a Lipschitz continuous boundary, one knows that functions in $\widetilde{H}^s(\Omega)$ coincide with $H^s(\R^n)$ functions vanishing a.e. in $\Omega^c$. Thus, we have $\psi-\Phi\in \widetilde{W}(0,T;\widetilde{H}^s(\Omega))$ and thus $u-\Phi\in \widetilde{W}(0,T;\widetilde{H}^s(\Omega))$. On the other hand, the condition \eqref{eq: equivalent characterization of ext cond} clearly implies $u=\Phi$ in $(\Omega_e)_T$.
        
        By the regularity assumptions of the involved functions and the above equivalent reformulation of the exterior condition, this means nothing else than that $v=u-\Phi\in \widetilde{W}(0,T;\widetilde{H}^s(\Omega))$ solves \eqref{eq: well-posedness linear case}, where the right hand side is given by
        \begin{equation}
        \label{eq: RHS in well-posedness for linear with nonzero ext cond}
            h-(\partial_t^2\Phi+(-\Delta)^s\partial_t \Phi+(-\Delta)^s\Phi+q\Phi)\in L^2(0,T;H^{-s}(\Omega))
        \end{equation}
        and the initial conditions by $u_0-\Phi(0)$, $u_1-\partial_t\Phi(0)$. To see the regularity condition in \eqref{eq: RHS in well-posedness for linear with nonzero ext cond} recall the estimate \eqref{eq: estimate a1} from the proof of Theorem~\ref{thm: well-posedness linear}. As this problem is well-posed the same holds for problem \eqref{eq: well-posedness linear case nonzero exterior cond}.\\

        \noindent \ref{item: sol indep of extension}: One easily sees that $u_1-u_2$ belongs to $\widetilde{W}(0,T;\widetilde{H}^s(\Omega))$ and is the unique solution of
        \begin{equation}
			\begin{cases}
				\LC \partial_t^2 +(-\Delta)^s\partial_t +(-\Delta)^s+q\RC u=0 &\text{ in }\Omega_T\\
				u=0     &\text{ in }(\Omega_e)_T,\\
				u(0)=0,\quad \partial_t u(0)=0  &\text{ in }\Omega.
			\end{cases}
		\end{equation}
        From Theorem~\ref{thm: well-posedness linear} we deduce $u_1=u_2$. Hence, we can conclude the proof.

    \end{proof}

    \subsection{Viscous wave equation with nonlinear perturbations}
    \label{subsec: well-posedness nonlinear}

    In this section we establish the well-posedness of the viscous wave equation with nonlinear perturbations. To this end we first discuss in Section~ \ref{subsec: diff of nemytskii} the continuity and differentiablity of the Nemytskii operator $f(u)$ under certain assumptions on $f$ (see Assumption~\ref{main assumptions on nonlinearities}). Then in Section~\ref{subsec: diff solution map} we invoke the implicit function theorem to get the desired well-posedness result and see that the solution map $S(\varphi)$ depends differentiabily on the exterior condition for small $\varphi$.
        
    \subsubsection{Differentiability of nonlinear perturbations}
    \label{subsec: diff of nemytskii}

    Next, we move on to the nonlinear problem. For this purpose let us specify a class of nonlinearities $f$ containing the one considered in the PDE \eqref{eq: nonlocal viscous wave eq nonlin pertur}. We start by recalling the notion of a Carath\'eodory function.
	
	\begin{definition}
		Let $U\subset \R^n$ be an open set. We say that  $f\colon U\times\R\to\R$ is a \emph{Carath\'odory function}, if it has the following properties:
		\begin{enumerate}[(i)]
			\item $\tau\mapsto f(x,\tau)$ is continuous for a.e. $x\in U$,
			\item $x\mapsto f(x,\tau)$ is measurable for all $\tau \in\R$.
		\end{enumerate}
	\end{definition}
	
	\begin{assumption}\label{main assumptions on nonlinearities}
		Let $f\colon \Omega\times \R\to\R$ a Carath\'eodory function satisfying the following conditions:
		\begin{enumerate}[(i)]
			\item\label{prop f} $f$ has partial derivative $\partial_{\tau}f$, which is a Carath\'eodory function,
            \item and there exists $a\in L^p(\Omega)$ such that
			\begin{equation}
				\label{eq: bound on derivative}
				\left|\partial_\tau f(x,\tau)\right|\lesssim a(x)+|\tau|^r
			\end{equation}
			for all $\tau\in\R$ and a.e. $x\in\Omega$. Here the exponents $p$ and $r$ satisfy the restrictions
			\begin{equation}\label{eq: cond on p}
				\begin{cases}
					n/s\leq p\leq \infty, &\, \text{if }\, 2s< n,\\
					2<p\leq \infty,  &\, \text{if }\, 2s= n,\\
					2\leq p\leq \infty, &\, \text{if }\, 2s > n,
				\end{cases}
			\end{equation} 
			and
			\begin{equation}
				\label{eq: cond on r}
				\begin{cases}
					0\leq r<\infty, &\, \text{if }\, 2s \geq n,\\
					0\leq r\leq \frac{2s}{n-2s}, &\, \text{if }\, 2s< n,
				\end{cases}
			\end{equation}
			respectively.
			
		\end{enumerate}
	\end{assumption}
	
	\begin{remark}
		An example of a nonlinearity \( f \) that satisfies the conditions in Assumption~\ref{main assumptions on nonlinearities} is the fractional power–type nonlinearity
\[
    f(x,\tau) = q(x)\,|\tau|^{r}\tau,
\]
where \( r \ge 0 \) satisfies \eqref{eq: cond on r} and \( q \in L^{\infty}(\Omega) \).  
The regularity conditions are clearly fulfilled. Moreover, a direct computation shows that
\begin{equation}\label{eq: f-power 1}
    \partial_\tau f(x,\tau) = (r+1)\, q(x)\, |\tau|^{r},
\end{equation}
and therefore, in this case, we may take \( a = 0 \) in \eqref{eq: bound on derivative}.
	\end{remark}
	
    Next, we state to auxiliary lemmas on the continuity and differentiability of Nemytskii operators.
	\begin{lemma}[Continuity of Nemytskii operators]
		\label{lemma: Continuity of Nemytskii operators}
		Let $\Omega\Subset \R^n$, $T>0$, $1\leq q,p<\infty$ and assume that $f\colon\Omega\times \R\to\R$ is a Carath\'eodory function satisfying
		\begin{equation}
			\label{eq: cond continuity nemytskii}
			|f(x,\tau)|\leq a+b|\tau|^{\alpha}
		\end{equation}
		for some constants $a,b\geq 0$ and $0<\alpha \leq \min(p,q)$. Then the Nemytskii operator $f$, defined by
		\begin{equation}
			\label{eq: Nemytskii operator}
			f(u)(x,t)\vcentcolon = f(x,u(x,t))
		\end{equation}
		for all measurable functions $u\colon \Omega_T\to\R$, maps continuously $L^q(0,T;L^p(\Omega))$ into $L^{q/\alpha}(0,T;L^{p/\alpha}(\Omega))$.
	\end{lemma}
	\begin{proof}
		The measurability and that the Nemytskii operator $f$ is well-defined, are immediate. Thus, we only need to check that it is continuous.
		
		Let $(u_n)_{n\in\N}\subset L^qL^p$ such that $u_n\to u$ in $L^qL^p$ as $n\to\infty$. By the converse of the dominated convergence theorem, we know that there is a subsequence, still denoted by $(u_n)$, and a function $g\in L^q((0,T))$ such that
		\begin{enumerate}[(a)]
			\item\label{cond 1} $u_n(t)\to u(t)$ in $L^p(\Omega)$ for a.e. $t$,
			\item\label{cond 2} $\|u_n(t)\|_{L^p(\Omega)}\leq g(t)$ for a.e. $t$.
		\end{enumerate}
		By \cite[Theorem 2.2]{ambrosetti1995primer} we know that $f$ is continuous from $L^p(\Omega)$ to $L^{p/\alpha}(\Omega)$ and so taking into account the estimate \eqref{eq: cond continuity nemytskii} as well as \ref{cond 1}, \ref{cond 2}, we deduce that
		\begin{enumerate}[(A)]
			\item $f(u_n(t))\to f(u(t))$ in $L^{p/\alpha}(\Omega)$ for a.e. $t$,
			\item $\|f(u_n(t))\|_{L^{p/\alpha}(\Omega)}\leq C(1+\|u_n(t)\|_{L^p(\Omega)}^{\alpha})\leq C(1+g(t)^{\alpha})\in L^{q/\alpha}((0,T))$.
		\end{enumerate}
		Applying the dominated convergence theorem shows that $f(u_n)\to f(u)$ in $L^{q/\alpha}L^{p/\alpha}$ as $n\to\infty$. Since this holds for every subsequence of $(u_n)_{n\in\N}$, we see that the whole sequence $(f(u_n))_{n\in\N}$ converges to $f(u)$. Hence, we can conclude the proof.
	\end{proof}
	
	\begin{lemma}[Differentiability of Nemytskii operators]
		\label{lemma: Differentiability of Nemytskii operators}
		Let $\Omega\Subset \R^n$ be a domain, $T>0$, $2<p<\infty$, $p-1\leq q<\infty$ and assume that $f\colon\Omega\times \R\to\R$ is a Carath\'eodory function with $f(\cdot,0)\in L^{\infty}(\Omega)$. Moreover, suppose that $f$ has partial derivative $\partial_\tau f$, which is a Carath\'eodory function and satisfies the estimate
		\begin{equation}
			\label{eq: derivative cond}
			\left|\partial_{\tau}f(x,\tau)\right|\leq a+b|\tau|^{p-2}
		\end{equation}
		for some constants $a,b\geq 0$. Then the Nemytskii operator $f$ is Fr\'echet differentiable as a map from $L^q(0,T;L^p(\Omega))$ to $L^{\frac{q}{p-1}}(0,T;L^{\frac{p}{p-1}}(\Omega))$ with differential
		\begin{equation}
			\label{eq: differential Nemytskii operator}
			df(u)h=\partial_\tau f(u)h.
		\end{equation}
	\end{lemma}
	\begin{proof}
		An integration of \eqref{eq: derivative cond} shows that $f$ satisfies the growth condition
		\begin{equation}
			\label{eq: growth of f}
			|f(x,\tau)|\leq c+d|\tau|^{p-1}
		\end{equation}
		for some constants $c,d\geq 0$. By Lemma~\ref{lemma: Continuity of Nemytskii operators} it follows that $f$ and $\partial_{\tau}f$ are continuous as mappings
		\begin{equation}
			\label{eq: cond f and partial f}
			\begin{split}
				f\colon L^q(0,T;L^p(\Omega))&\to L^{\frac{q}{p-1}}(0,T;L^{\frac{p}{p-1}}(\Omega))\\
				\partial_{\tau}f\colon L^q(0,T;L^p(\Omega))&\to L^{\frac{q}{p-2}}(0,T;L^{\frac{p}{p-2}}(\Omega)).
			\end{split}
		\end{equation}
		Next, let $u,h\in L^qL^p$ and define 
		\[
		\omega(u,h)=\|f(u+h)-f(u)-\partial_\tau f(u) h\|_{L^{\frac{q}{p-1}}L^{\frac{p}{p-1}}}.
		\]
		By \cite[eq. (2.8)]{ambrosetti1995primer}, we know
		\[
		\omega(u,h)\leq \left\|\|h\|_{L^p}\left\|\int_0^1\left(\partial_\tau f(u+\xi h)-\partial_\tau f(u)\right)\,d\xi\right\|_{L^{\frac{p}{p-2}}}\right\|_{L^{\frac{q}{p-1}}}.
		\]
		Observing that 
		\[
		\frac{p-1}{q}=\frac{1}{q}+\frac{p-2}{q},
		\]
		we get by H\"older's and Minkowski's inequality
		\[
		\begin{split}
			\omega(u,h)&\leq \|h\|_{L^qL^p}\left\|\int_0^1\left(\partial_\tau f(u+\xi h)-\partial_\tau f(u)\right)d\xi\right\|_{L^{\frac{q}{p-2}}L^{\frac{p}{p-2}}}\\
			&\leq \|h\|_{L^qL^p}\int_0^1\left\|\partial_\tau f(u+\xi h)-\partial_\tau f(u)\right\|_{L^{\frac{q}{p-2}}L^{\frac{p}{p-2}}}d\xi.
		\end{split}
		\]
		By \eqref{eq: cond f and partial f} the second factor goes to zero as $h\to 0$ in $L^qL^p$ and therefore we get $\omega(u,h)=o(\|h\|_{L^qL^p})$. Hence, $f$ is a differentiable map from $L^qL^p$ to $L^{\frac{q}{p-1}}L^{\frac{p}{p-1}}$.
	\end{proof}

    \subsubsection{Differentiability of solution map to the nonlinear problem}
    \label{subsec: diff solution map}

    Now with the tools from the preceding section at our disposal, we can show that the viscous wave equation with nonlinear perturbations is well-posed for small exterior conditions.
	\begin{theorem}
		\label{prop: differentiability of solution map}
		Let $\Omega\subset\R^n$ be a bounded Lipschitz domain, $T>0$ and $s>0$. Suppose that $f$ satisfies $f(0)=0$ and Assumption~\ref{main assumptions on nonlinearities} with $r>0$ and $a\in L^{\infty}(\Omega)$. Then there exist neighborhoods $U_0\subset \widetilde{W}^s_\text{rest}((\Omega_e)_T)$, $U_1\subset \widetilde{W}_{ext}(0,T;\widetilde{H}^s(\Omega))$ of the origin with the property that, for any $\varphi\in U_0$, the problem
		\begin{equation}
        \label{eq: nonlinear pde with nonzero ext cond}
		\begin{cases}
			\partial_t^2u +(-\Delta)^s\partial_t u+(-\Delta)^su+f(u)= 0&\text{ in }\Omega_T\\
			u=\varphi&\text{ in }(\Omega_e)_T,\\
			u(0)=0,\quad \partial_t u(0)=0  &\text{ in }\Omega
		\end{cases}
		\end{equation}
		has a unique solution $u\in U_1$ in the sense that:
        \begin{enumerate}[(i)]
            \item\label{PDE nonlinear} for all $v\in\widetilde{H}^s(\Omega)$ one has
            \begin{equation}
            \label{eq: weak solution linear nonzero ext cond}
            \begin{split}
                 &\frac{d}{dt}\,\langle \partial_t u,v\rangle_{L^2(\Omega)}+\langle (-\Delta)^{s/2}\partial_t u,(-\Delta)^{s/2}v\rangle_{L^2(\R^n)}\\
                 &+\langle (-\Delta)^{s/2}u,(-\Delta)^{s/2}v\rangle_{L^2(\R^n)}+\langle f(u),v\rangle_{L^2(\Omega)}=0
            \end{split}
        \end{equation}
         in the sense of $\distr((0,T))$;
            \item\label{exterior cond nonlinear} $u=\varphi$ in $(\Omega_e)_T$;
            \item\label{initial values nonlinear} $u(0)=0$ in $\widetilde{H}^s(\Omega)$ and $\partial_t u(0)=0$ in $L^2(\Omega)$.
        \end{enumerate}
        Moreover, the map $U_0\ni\varphi\mapsto S(\varphi)\in U_1$, which assigns to each exterior condition $\varphi\in U_0$ the corresponding unique solution $S(\varphi)\vcentcolon = u$ of \eqref{eq: nonlinear pde with nonzero ext cond}, is $C^1$ in the Fr\'echet sense.
	\end{theorem}
    \begin{remark}
        Heuristically, the map $S$ assigns to each "small" exterior condition $\varphi$ its unique "small" solution $S(\varphi)=u$ and we refer to $S$ as the \emph{solution map} associated with the problem  \eqref{eq: nonlinear pde with nonzero ext cond}. 
    \end{remark}

	\begin{proof}
		
		Let us start by defining the following Banach spaces
		\begin{equation}
			\label{eq: function spaces differentiability}
			\begin{split}
				E_0&=\widetilde{W}^s_\text{rest}((\Omega_e)_T),\quad E_1=\widetilde{W}_{ext}(0,T;\widetilde{H}^s(\Omega)),\\
				V_0&=\widetilde{H}^s(\Omega),\quad V_1=L^2(\Omega),\quad V_2=\widetilde{W}^s_\text{rest}((\Omega_e)_T),\quad  V_3=L^{2}(0,T;H^{-s}(\Omega)),
			\end{split}
		\end{equation}
		and define the map $F: E_0\times E_1\to \prod_{j=0}^3 V_j$ via
		\begin{equation}
			\label{eq: F for diff}
			F(\varphi,u) =\LC u(0), \p_tu(0),(u-\varphi)|_{(\Omega_e)_T},(\p_t^2+(-\Delta)^s\partial_t+(-\Delta)^s)u+f(u)\RC .
		\end{equation}		
        We now wish to argue that $F$ is well-defined. For this purpose we first show that there holds
        \begin{equation}
			\label{eq: estimate Hs}
			\begin{split}
				\|f(\psi)\|_{L^2(\Omega_T)}
				\lesssim\, &(\|a\|_{L^{\infty}(\Omega)}\|\psi\|_{L^{\infty}(0,T;H^s(\R^n))} +\|\psi\|^{r+1}_{L^{\infty}(0,T;H^s(\R^n))}).
			\end{split}
		\end{equation}
        for any $\psi\in C([0,T];H^s(\R^n))$. First of all by the fundamental theorem of calculus, Assumption~\ref{main assumptions on nonlinearities} and $f(0)=0$ we have
		\begin{equation}
			\label{eq:estimate on f}
			|f(x,s)|\leq \left|\int_0^s \partial_{\tau}f(x,\tau)\,d\tau\right|\leq C\LC a(x)|s|+|s|^{r+1}\RC,
		\end{equation}
		where $a \in L^{\infty}(\Omega)$ is nonnegative. This ensures that we have
		\[
		|f(\psi(t))|\leq C\LC a|\psi(t)|+|\psi(t)|^{r+1}\RC.
		\]
		Hence, we get
		\begin{equation}
			\label{eq: L2 estimate f}
			\|f(\psi(t))\|_{L^2(\Omega)}\leq C\LC \|a\psi(t)\|_{L^2(\Omega)}+\|\psi(t)\|_{L^{2(r+1)}(\Omega)}^{r+1}\RC,
		\end{equation}
		for $0\leq t\leq T$. By the boundedness of $a$, we obtain
		\begin{equation}
			\label{eq: linear growth term}
			\|a\psi(t)\|_{L^2(\Omega)}\leq \|a\|_{L^{\infty}(\Omega)}\|\psi(t)\|_{H^s(\R^n)}.
		\end{equation}
		Note that the conditions on the exponent $r$ yield
		\[
		\begin{cases}
			1\leq 1+r<\infty, &\, \text{if }\, 2s\geq n,\\
			1\leq 1+r\leq \frac{n}{n-2s} &\, \text{if }\, 2s< n.
		\end{cases}
		\]
		If $2s> n$, then the Sobolev embedding $H^s(\R^n)\hookrightarrow L^{\infty}(\R^n)$ yields
		\begin{equation}
			\label{eq: supercritical case}
			\|\psi(t)\|_{L^{2(r+1)}(\Omega)}^{r+1}\leq C\|\psi(t)\|_{H^s(\R^n)}^{r+1}.
		\end{equation}
		In the critical case $2s=n$, we can apply \cite{Ozawa} to obtain
		\begin{equation}
			\label{eq: critical case}
			\|\psi(t)\|_{L^{2(r+1)}(\Omega)}^{r+1}\leq C\|(-\Delta)^{s/2}\psi(t)\|_{L^2(\R^n)}^r\|\psi(t)\|_{L^2(\R^n)}.
		\end{equation}
		In the subcritical case $2s<n$, we apply the Hardy--Littlewood--Sobolev lemma to deduce
		\begin{equation}
			\label{eq: subcritical case}
			\|\psi(t)\|_{L^{2(r+1)}(\Omega)}^{r+1}\leq C\|\psi(t)\|_{L^{\frac{2n}{n-2s}}(\Omega)}^{r+1}\leq C \|(-\Delta)^{s/2}\psi(t) \|_{L^2(\R^n)}^{r+1}.
		\end{equation}
		As $\psi\in C([0,T];H^s(\R^n))$, we get by the continuity of the fractional Laplacian the desired estimate \eqref{eq: estimate Hs}.
  
		Now we can show that $F$ is well-defined. Since $u\in \widetilde{W}_{ext}(0,T;\widetilde{H}^s(\Omega))$, we have by definition $u(0)\in\widetilde{H}^s(\Omega)$ and $\partial_t u(0)\in L^2(\Omega)$. Thus, the first two entries of $F$ are well-defined. On the other hand as $u=v+\psi\in \widetilde{W}_{ext}(0,T;\widetilde{H}^s(\Omega))$, we have $(u-\varphi)|_{(\Omega_e)_T}=(\psi-\varphi)|_{(\Omega_e)_T}\in \widetilde{W}^s_\text{rest}((\Omega_e)_T)$. Finally, using the mapping properties of $(-\Delta)^s$, Lemma~\ref{lemma: embeddings} and \eqref{eq: estimate Hs}, one easily sees that $\partial_t^2 u+(-\Delta)^s\partial_t+(-\Delta)^s u+f(u)\in L^2(0,T;H^{-s}(\Omega))$ for $u\in\widetilde{W}_{ext}(0,T;\widetilde{H}^s(\Omega))$.

		We next show that $F$ is (Fr\'echet) differentiable. Note that all appearing operators, up to $f(u)$, are linear bounded operators and hence differentiable. Thus, it remains to show that $f(u)$ is differentiable from $E_1 \to V_3$. 
		
		\medskip

		\noindent \textit{Case $2s<n$.} By Assumption~\ref{main assumptions on nonlinearities}, $a\in L^{\infty}(\Omega)$ and $r>0$, we see that all conditions in Lemma~\ref{lemma: Differentiability of Nemytskii operators} are satisfied, when $p=r+2$ and $1\leq q<\infty$ is any number satisfying $r+1\leq q<\infty$. Hence, the Nemytskii operator $f(u)$ is differentiable as a map from $L^q(0,T;L^{r+2}(\Omega))$ to $L^{\frac{q}{r+1}}(0,T;L^{\frac{r+2}{r+1}}(\Omega))$ with differential
		\begin{equation}
			\label{eq: derivative of nonlinearity}
			df(u)h=\partial_{\tau}f(u)h.
		\end{equation}
		Next, recall that $r$ satisfies the condition \eqref{eq: cond on r} so that 
		\[
		1< 1+r\leq 1+\frac{2s}{n-2s}=\frac{n}{n-2s}.
		\]
		This implies
		\[
		2<2+r\leq 2+\frac{2s}{n-2s}=\frac{2n-2s}{n-2s}<\frac{2n}{n-2s}.
		\]
		Thus, by the Sobolev embedding and boundedness of $\Omega$ we get
		\begin{equation}
			\label{eq: Hs to L r2}
			H^s(\R^n)\hookrightarrow L^{\frac{2n}{n-2s}}(\Omega)\hookrightarrow L^{r+2}(\Omega).
		\end{equation}
		On the other hand the conjugate exponent
		\[
		(r+2)'=\frac{r+2}{r+1}=1+\frac{1}{r+1}
		\]
		fulfills
		\[
		2>(r+2)'\geq 1+\frac{n-2s}{n}=\frac{2(n-s)}{n}.
		\]
		Next, observe that
		\[
		\frac{2n}{n+2s}< \frac{2(n-s)}{n} \quad \Leftrightarrow \quad n^2< n^2+sn-2s^2
		\]
		and thus by the Sobolev embedding $L^{\frac{2n}{n+2s}}(\R^n)\hookrightarrow H^{-s}(\Omega)$ we obtain
		\begin{equation}
			\label{eq: right embed}
			L^{(r+2)'}(\Omega)\hookrightarrow H^{-s}(\Omega).
		\end{equation}
		Combining \eqref{eq: Hs to L r2} and \eqref{eq: right embed}, we get that $u\mapsto f(u)$ is differentiable as a map from $L^q(0,T;H^s(\R^n))$ to $L^{\frac{q}{r+1}}(0,T;H^{-s}(\Omega))$. Choosing $q=2(r+1)\geq 2$ and using the embedding $\widetilde{W}(0,T;\widetilde{H}^s(U))\hookrightarrow C([0,T];\widetilde{H}^s(U))$ for any open set $U\subset\R^n$, we see that $f(u)$ is differentiable as a map from $E_1$ to $V_3$. Next, we assert that the differential, given by \eqref{eq: derivative of nonlinearity}, is continuous as a map from $E_1$ to $L(E_1,V_3)$. This then establishes that $f$ is $C^1$ as a map from $E_1$ to $V_3$. By Lemma~\ref{lemma: Continuity of Nemytskii operators} with $\alpha=r$, we know that $\partial_\tau f (u)$ is continuous as a map from $L^q(0,T;L^{r+2}(\Omega))$ to $L^{\frac{q}{r}}(0,T;L^{\frac{r+2}{r}}(\Omega))$, when $q\geq r$. Now, let us choose $q$ such that $q\geq 2\max(r,1)$ and observe that by Assumption~\ref{main assumptions on nonlinearities} there holds
		\[
		\frac{r+2}{r}=1+\frac{2}{r}\in \left[n/s-1,\infty\right).
		\]
        Therefore, we can define
        \begin{equation}
        \label{eq: def of p}
            \frac{1}{p}=\frac{n+2s}{2n}-\frac{r}{r+2}>0.
        \end{equation}
        One may observe that 
        \[
            \frac{n+2s}{2n}-\frac{s}{n-s}\geq \frac{n-2s}{2n}
        \]
       and hence one has
        \begin{equation}
        \label{eq: range of p}
            \frac{1}{p}\geq \frac{n+2s}{2n}-\frac{s}{n-s}\geq \frac{n-2s}{2n}.
        \end{equation}
		Now let $u_k\in E_1$, $k\in\N$, converge to some $u\in E_1$ and fix some functions $v\in E_1$, $\psi\in L^2(0,T;\widetilde{H}^s(\Omega))$. Then by H\"olders and Sobolev's inequality we can estimate
        \[
            \begin{split}
                &|\langle (\partial_\tau f(u_k)-\partial_\tau f(u))v,\psi\rangle|\\
                &\lesssim \|(\partial_\tau f(u_k)-\partial_\tau f(u))v\|_{L^2(0,T;L^{\frac{2n}{n+2s}}(\Omega))}\|\psi\|_{L^2(0,T;L^{\frac{2n}{n-2s}}(\Omega))}\\
                &\overset{\eqref{eq: def of p}}{\lesssim} \|\partial_\tau f(u_k)-\partial_\tau f(u)\|_{L^2(0,T;L^{\frac{r+2}{r}}(\Omega))}\|v\|_{L^{\infty}(0,T;L^p(\Omega))}\|\psi\|_{L^2(0,T;\widetilde{H}^s(\Omega))}\\
                &\overset{\eqref{eq: range of p}}{\lesssim} \|\partial_\tau f(u_k)-\partial_\tau f(u)\|_{L^2(0,T;L^{\frac{r+2}{r}}(\Omega))}\|v\|_{L^{\infty}(0,T;L^{\frac{2n}{n-2s}}(\Omega))}\|\psi\|_{L^2(0,T;\widetilde{H}^s(\Omega))}\\
                &\lesssim \|\partial_\tau f(u_k)-\partial_\tau f(u)\|_{L^2(0,T;L^{\frac{r+2}{r}}(\Omega))}\|v\|_{L^{\infty}(0,T;H^s(\R^n))}\|\psi\|_{L^2(0,T;\widetilde{H}^s(\Omega))}\\
                &\lesssim \|\partial_\tau f(u_k)-\partial_\tau f(u)\|_{L^2(0,T;L^{\frac{r+2}{r}}(\Omega))}\|v\|_{E_1}\|\psi\|_{L^2(0,T;\widetilde{H}^s(\Omega))}\\
                &\lesssim \|\partial_\tau f(u_k)-\partial_\tau f(u)\|_{L^{\frac{q}{r}}(0,T;L^{\frac{r+2}{r}}(\Omega))}\|v\|_{E_1}\|\psi\|_{L^2(0,T;\widetilde{H}^s(\Omega))}
            \end{split}
         \]
		This shows that
        \begin{equation}
        \label{eq: continuity of derivative}
         \|\partial_\tau f(u_k)-\partial_\tau f(u)\|_{L(E_1,V_3)}\lesssim \|(\partial_\tau f(u_k)-\partial_\tau f(u))\|_{L^{\frac{q}{r}}(0,T;L^{\frac{r+2}{r}}(\Omega))}.
        \end{equation}
        Now, by the continuity of $\partial_\tau f(u)$ from $L^q(0,T;L^{r+2}(\Omega))$ to $L^{\frac{q}{r}}(0,T;L^{\frac{r+2}{r}}(\Omega))$ and the embedding $E_1\hookrightarrow L^{\Bar{p}}(0,T;H^s(\R^n))$ for any $1\leq \Bar{p}\leq \infty$, we see that \eqref{eq: continuity of derivative} goes to zero as $k\to\infty$ and hence the differential is continuous as we wanted to show.
		\medskip
		
		\noindent  \textit{Case $2s\geq n$.} After recalling that in this case we have $H^s(\R^n)\hookrightarrow L^p(\Omega)$ for any $2\leq p<\infty$ and $L^q(\Omega)\hookrightarrow H^{-s}(\Omega)$ for any $1<q\leq 2$ (see \cite{Ozawa} for supercritical Sobolev embedding), one can argue similarly as in the subcritical case $2s<n$.\\
		
		\noindent Hence, $F$ is a $C^1$ map. Next note that $F(0,0)=0$ . Now, the derivative of $F$ at the origin in the $u$-variable is
		\[
		\partial_u F(0,0)v = (v(0),\partial_t v(0),v|_{(\Omega_e)_T}, (\p_t^2+(-\Delta)^s\partial_t+(-\Delta)^s+\p_\tau f(0))v)
		\]
		for $v\in E_1$. This map is a linear, bounded and invertible operator from $E_1\to \prod_{j=0}^3 V_j$. To see this, consider the problem
		\begin{equation}
        \label{eq: linearized equation frechet}
		\begin{cases}
			( \p_t^2 +(-\Delta)^s\partial_t+ (-\Delta)^s + \p_\tau f(0))v=h &\text{ in }\Omega_T\\
			v=\psi &\text{ in } (\Omega_e)_T,\\
			v=v_0,\quad \p_tv=v_1 &\text{ in }\Omega
		\end{cases}
		\end{equation}
    for $v_0\in\widetilde{H}^s(\Omega)$, $v_1\in L^2(\Omega)$, $\psi\in W^s_{rest}((\Omega_e)_T)$ and $h\in L^2(0,T;H^{-s}( \Omega))$. The well-posedness of this problem follows from Corollary~\ref{Well-posedness with nonzero exterior conditions}.
		
		Now, the implicit function theorem on Banach spaces \cite[Theorem~2.3]{ambrosetti1995primer} yields that there exist neighborhoods $U_0\subset E_0$, $U_1\subset E_1$ containing the origin and a map $S\in C^1(U_0,E_1)$ such that 
		\begin{enumerate}[(i)]
		    \item $F(\varphi, S(\varphi))=0$ for all $\varphi\in U_0$,
            \item $F(\varphi,u)=0$ for some $(\varphi,u)\in U_0\times U_1$, then $u=S(\varphi)$.
		\end{enumerate}
        One easily sees that $u=S(\varphi)$, for $\varphi\in U_0$, satisfies the conditions \ref{PDE nonlinear}--\ref{initial values nonlinear}.
	\end{proof}
	
    \section{Inverse problem for the linear viscous wave equation}
    \label{sec: solving linear IP}

    In this section we move on to the inverse problem for the viscous wave equation with linear perturbations. First in Section~ \ref{DN map and runge in linear problem} we introduce rigorously the corresponding DN map (Definition~\ref{DN map linear case}) and then prove the Runge approximation (Proposition~\ref{prop: runge}) in $L^2(0,T;\widetilde{H}^s(\Omega))$. Then in Section~\ref{subsec: unique determination linear} we present the proof of Theorem~\ref{Thm: main linear} after establishing a suitable integral identity in Lemma~\ref{lemma: integral identity}.
    
	\subsection{DN map and Runge approximation for the linear problem}
 \label{DN map and runge in linear problem}

    \begin{definition}
    \label{DN map linear case}
        Let $\Omega\subset\R^n$ be a bounded Lipschitz domain, $T>0$ and $s>0$. Suppose that the (real valued) function $q\in L^1_{loc}(\Omega_T)$ satisfies the conditions in Theorem~\ref{thm: well-posedness linear}. Then we define the \emph{Dirichlet to Neumann map} $\Lambda_q$ by
        \begin{equation}
        \label{eq: DN map linear case}
            \langle \Lambda_q\varphi,\psi\rangle=\int_{\R^n_T} (-\Delta)^{s/2} u (-\Delta)^{s/2} \psi\,dxdt+\int_{\R^n_T}(-\Delta)^{s/2}\partial_t u(-\Delta)^{s/2}\psi\,dxdt
        \end{equation}
        for all $\varphi,\psi\in C_c^{\infty}((\Omega_e)_T)$. Here, $u\in\widetilde{W}_{ext}(0,T;\widetilde{H}^s(\Omega))$ denotes the unique solution of 
        \[
        \begin{cases}
				\LC \partial_t^2 +(-\Delta)^s\partial_t +(-\Delta)^s+q\RC u=0 &\text{ in }\Omega_T\\
				u=\varphi     &\text{ in }(\Omega_e)_T,\\
				u(0)=0,\quad \partial_t u(0)=0  &\text{ in }\Omega.
			\end{cases}
        \]
    \end{definition}
 
	\begin{proposition}[Runge approximation]
		\label{prop: runge}
		Let $\Omega\subset\R^n$ be a bounded Lipschitz domain, $W\subset\Omega_e$ a given measurement set, $T>0$ and $s>0$ a non-integer. Suppose that the (real valued) function $q\in L^1_{loc}(\Omega_T)$ satisfies the conditions in Theorem~\ref{thm: well-posedness linear}. Consider the \emph{Runge set} 
		\begin{align*}
			\mathcal{R}_W:=\left\{ u_{\varphi}-\varphi\, : \, \varphi\in C^\infty_c (W_T) \right\},
		\end{align*}
		where $u_{\varphi}\in \widetilde{W}_{ext}(0,T;\widetilde{H}^s(\Omega))$ is the unique solution to 
		\[
		\begin{cases}
			(\partial_t^2 +(-\Delta)^s\partial_t+(-\Delta)^s + q)u = 0 &\text{ in }\Omega_T \\
			u=\varphi  &\text{ in }(\Omega_e)_T,\\
			u(0)=0, \quad \partial_t u(0)=0 &\text{ in }\Omega.
		\end{cases}
		\] 
		Then $\mathcal{R}_W$ is dense in $L^2(0,T;\widetilde{H}^s(\Omega))$.
	\end{proposition}
	
	\begin{proof}
		Since $\mathcal{R}_W\subset L^2(0,T;\widetilde{H}^s(\Omega))$ is a subspace it is enough by the Hahn--Banach theorem to show that if $F\in L^2(0,T;H^{-s}(\Omega))$ vanishes on $\mathcal{R}_W$, then $F=0$. Hence, choose any $F\in L^2(0,T;H^{-s}(\Omega))$ and assume that 
		\begin{equation}
			\langle F,u_{\varphi}-\varphi\rangle=0\quad\text{for all}\quad \varphi\in C_c^{\infty}(W_T).
		\end{equation}
		Next, let $w_F\in\widetilde{W}(0,T;\widetilde{H}^s(\Omega))$ be the unique solution to the adjoint equation 
		\begin{equation}
			\begin{cases}
				(\partial_t^2 -(-\Delta)^s\partial_t+(-\Delta)^s+q)w = F&\text{ in }\Omega \times (0,T)\\
				w=0  &\text{ in }\Omega_e\times (0,T),\\
				w(T)=0, \quad \partial_t w(T)=0 &\text{ in }\Omega,
			\end{cases}
		\end{equation}
		which exists by Theorem~\ref{thm: well-posedness linear} with $q$ replaced by $q^\star$ in \eqref{eq: well-posedness linear case} and a subsequent time reversal of the solution. Next, let us note that the integration by parts formula gives us
        \begin{equation}
        \label{eq: integration by parts 2nd order der}
            \begin{split}
            \int_0^T \langle \partial_t^2 v,w\rangle\,dt&=\int_0^T\langle \partial_t^2 w,v\rangle\,dt +\langle \partial_t v(T),w(T)\rangle -\langle \partial_t w(T),v(t)\rangle\\
            &\quad -(\langle \partial_t v(0),w(0)\rangle-\langle \partial_t w(0),v(0)\rangle)
            \end{split}
        \end{equation}
        and
        \begin{equation}
        \label{eq: integration by parts 1st order der}
            \begin{split}
                &\int_0^T \langle (-\Delta)^{s/2}\partial_t v,(-\Delta)^{s/2}w\rangle\,dt=-\int_0^T \langle (-\Delta)^{s/2}v,(-\Delta)^{s/2}\partial_t w\rangle\,dt\\
            &\quad +\langle (-\Delta)^{s/2}v(T),(-\Delta)^{s/2}w(T)\rangle-\langle (-\Delta)^{s/2}v(0),(-\Delta)^{s/2}w(0)\rangle
            \end{split}
        \end{equation}
        for all $v,w\in\widetilde{W}(0,T;\widetilde{H}^s(\Omega))$. 
        By a density argument, the PDEs for $u-\varphi$ and $w_F$ hold in the $L^2(0,T;H^{-s}(\Omega))$ sense, with the help of \eqref{eq: integration by parts 2nd order der}, \eqref{eq: integration by parts 1st order der} and the vanishing initial and terminal conditions, respectively, we may compute
		\[
		\begin{split}
			0&=\left\langle F,u_{\varphi}-\varphi\right\rangle \\
			&=\left\langle (\partial_t^2 -(-\Delta)^s\partial_t+(-\Delta)^s +q)w_F,u_{\varphi}-\varphi\right\rangle\\
			&= \left\langle (\partial_t^2 +(-\Delta)^s\partial_t+(-\Delta)^s +q)(u_{\varphi}-\varphi),w_F\right\rangle\\
            &= -\left\langle (\partial_t^2 +(-\Delta)^s\partial_t+(-\Delta)^s +q)\varphi,w_F\right\rangle\\
            &=-\langle (-\Delta)^s\partial_t\varphi+(-\Delta)^s\varphi,w_F\rangle\\
			&=\left\langle(-\Delta)^s(w_F-\partial_t w_F),\varphi \right\rangle,
		\end{split}
		\]
		for all $\varphi\in C_c^{\infty}(W_T)$. This implies that $\widetilde{w}_F=w_F-\partial_tw_F\in L^2(0,T;\widetilde{H}^s(\Omega))$ satisfies
		\[
		(-\Delta)^s \widetilde{w}_F=\widetilde{w}_F=0\quad\text{in}\quad W_T.
		\]
		By the unique continuation property of the fractional Laplacian \cite[Theorem~1.2]{GSU20}, this gives $\widetilde{w}_F=0$ in $\R^n_T$. By construction we have $w_F\in H^1(0,T;\widetilde{H}^s(\Omega))$ and hence $w_F(x,\cdot)\in H^1((0,T))$ for a.e. $x\in \Omega$. Then as $\widetilde{w}_F=0$ in $\R^n_T$ we know that $w_F(x,\cdot)$ solves
        \begin{equation}
        \label{eq: ODE}
            \begin{cases}
                \partial_t w_F=w_F,\\
                w_F(T)=0
            \end{cases}
        \end{equation}
        and hence we may conclude that $w_F(x,\cdot)=0$ for a.e. $x\in\Omega$. Thus, we deduce that $w_F=0$ and therefore it follows that $F=0$ as we wanted to show.
	\end{proof}

    \subsection{Unique determination of linear perturbations}
    \label{subsec: unique determination linear}

    In this section we give the proof of Theorem~\ref{Thm: main linear}. We first deduce a suitable integral identity, which plays the role of the Alessandrini identity in the elliptic case.

    As already observed, below we will make use of the following simple fact: The function $u\in\widetilde{W}_{ext}(0,T;\widetilde{H}^s(\Omega))$ is the unique solution of 
    \begin{equation}
    \begin{cases}
			(\partial_t^2 +(-\Delta)^s\partial_t+(-\Delta)^s + q)u = 0 &\text{ in }\Omega_T \\
			u=\varphi  &\text{ in }(\Omega_e)_T,\\
			u(0)=0, \quad \partial_t u(0)=0 &\text{ in }\Omega,
	\end{cases}
    \end{equation}
    if and only if $u^\star \in \widetilde{W}_{ext}^{\star}(0,T;\widetilde{H}^s(\Omega))$ is the unique solution of
    \[
    \begin{cases}
			(\partial_t^2 -(-\Delta)^s\partial_t+(-\Delta)^s + q^\star)v = 0 &\text{ in }\Omega_T \\
			v=\varphi^\star  &\text{ in }(\Omega_e)_T,\\
			v(T)=0, \quad \partial_t v(T)=0 &\text{ in }\Omega.
	\end{cases}
    \]
    Here, $\widetilde{W}^\star(0,T;\widetilde{H}^s(\Omega))$ denotes the space $\widetilde{W}(0,T;\widetilde{H}^s(\Omega))+\widetilde{W}^{\star}_*(0,T;H^s(\R^n))$ with
    \[
        \widetilde{W}^{\star}_*(0,T;H^s(\R^n))=\{v\in\widetilde{W}(0,T;H^s(\R^n))\,;\,v(T)\in\widetilde{H}^s(\Omega),\,\partial_t v(T)\in L^2(\Omega))\}.
    \]
    \begin{lemma}[Integral identity]
    \label{lemma: integral identity}
        Let $\Omega\subset\R^n$ be a bounded Lipschitz domain, $T>0$ and $s>0$. Suppose that the (real valued) function $q\in L^1_{loc}(\Omega_T)$ satisfies the conditions in Theorem~\ref{thm: well-posedness linear}. Then there holds
        \begin{equation}
        \label{eq: Alessandrini identity}
            \langle(\Lambda_{q_1}-\Lambda_{q_2^{\star}})\varphi_1,\varphi_2^{\star}\rangle=\int_{\Omega_T}(q_1-q_2^{\star})(u_1-\varphi_1)(u_2-\varphi_2)^{\star}dxdt
        \end{equation}
        for all $\varphi_1,\varphi_2\in C_c^{\infty}((\Omega_e)_T)$, where $u_j\in\widetilde{W}_{ext}(0,T;\widetilde{H}^s(\Omega))$ is the unique solution of 
        \begin{equation}
        \label{eq: viscous wave for int id}
    \begin{cases}
			(\partial_t^2 +(-\Delta)^s\partial_t+(-\Delta)^s + q_j)u = 0 &\text{ in }\Omega_T \\
			u=\varphi_j  &\text{ in }(\Omega_e)_T,\\
			u(0)=0, \quad \partial_t u(0)=0 &\text{ in }\Omega
	\end{cases}
    \end{equation}
    for $j=1,2$.
    \end{lemma}

    \begin{proof}
        Let $Q_1,Q_2$ be two potentials satisfying the assumptions of Theorem~\ref{thm: well-posedness linear} and denote by $U_1,U_2^{\star}$ the unique solutions of 
        \begin{equation}
        \label{eq: for Q1}
            \begin{cases}
			(\partial_t^2 +(-\Delta)^s\partial_t+(-\Delta)^s + Q_1)u = 0 &\text{ in }\Omega_T \\
			u=\varphi_1  &\text{ in }(\Omega_e)_T,\\
			u(0)=0, \quad \partial_t u(0)=0 &\text{ in }\Omega
	\end{cases}
        \end{equation}
        and 
        \begin{equation}
        \label{eq: for Q2}
            \begin{cases}
			(\partial_t^2 -(-\Delta)^s\partial_t+(-\Delta)^s + Q_2)v = 0&\text{ in }\Omega_T \\
			v=\varphi_2^{\star}  &\text{ in }(\Omega_e)_T,\\
			v(T)=0, \quad \partial_t v(T)=0 &\text{ in }\Omega,
	\end{cases}
        \end{equation}
        respectively. Clearly, $U_2^\star$ is the time reversal of the solution $U_2$ solving
         \begin{equation}
        \label{eq: for Q2star}
            \begin{cases}
			(\partial_t^2 +(-\Delta)^s\partial_t+(-\Delta)^s + Q_2^\star)w = 0&\text{ in }\Omega_T \\
			w=\varphi_2  &\text{ in }(\Omega_e)_T,\\
			w(0)=0, \quad \partial_t w(0)=0 &\text{ in }\Omega.
	\end{cases}
        \end{equation}
        By \eqref{eq: integration by parts 2nd order der}, we know that there holds
        \begin{equation}
        \label{eq: 2nd order time derivatives integral id}
            \int_0^T \langle \partial_t^2(U_1-\varphi_1),(U_2-\varphi_2)^\star\rangle\,dt=\int_0^T \langle \partial_t^2(U_2-\varphi_2)^\star,U_1-\varphi_1\rangle\,dt.
        \end{equation}
        Therefore, we may compute
         \allowdisplaybreaks
         \small
        \begin{equation}
        \begin{split}
            &\int_{\Omega_T}(Q_1-Q_2)(U_1-\varphi_1)(U_2-\varphi_2)^\star\,dxdt\\
        &=-\int_0^T\langle(\partial_t^2+(-\Delta)^s\partial_t +(-\Delta)^s)(U_1-\varphi_1),(U_2-\varphi_2)^\star\rangle\,dt\\
        &+\int_0^T\langle(\partial_t^2-(-\Delta)^s\partial_t +(-\Delta)^s)(U_2-\varphi_2)^\star,U_1-\varphi_1\rangle\,dt\\
        &-\int_0^T\langle(\partial_t^2+(-\Delta)^s\partial_t +(-\Delta)^s)\varphi_1,(U_2-\varphi_2)^\star\rangle\,dt\\
        &+\int_0^T\langle(\partial_t^2-(-\Delta)^s\partial_t +(-\Delta)^s)\varphi_2^\star,U_1-\varphi_1\rangle\,dt\\
        &\overset{\eqref{eq: 2nd order time derivatives integral id}}{=}
        -\int_0^T\langle((-\Delta)^s\partial_t +(-\Delta)^s)(U_1-\varphi_1),(U_2-\varphi_2)^\star\rangle\,dt\\
        &+\int_0^T\langle(-(-\Delta)^s\partial_t +(-\Delta)^s)(U_2-\varphi_2)^\star,U_1-\varphi_1\rangle\,dt\\
        &-\int_0^T\langle((-\Delta)^s\partial_t +(-\Delta)^s)\varphi_1,(U_2-\varphi_2)^\star\rangle\,dt+\int_0^T\langle(-(-\Delta)^s\partial_t +(-\Delta)^s)\varphi_2^\star,U_1-\varphi_1\rangle\,dt.
        \end{split}
        \end{equation}
        \normalsize
        In the second equality sign we also used the support conditions of $\varphi_j$. As the first two terms compensate each other, using integration by parts we get
        \small
        \begin{equation}
        \begin{split}
        &\int_{\Omega_T}(Q_1-Q_2)(U_1-\varphi_1)(U_2-\varphi_2)^\star\,dxdt\\
        &= -\int_0^T\langle((-\Delta)^s\partial_t +(-\Delta)^s)\varphi_1,(U_2-\varphi_2)^\star\rangle\,dt+\int_0^T\langle(-(-\Delta)^s\partial_t +(-\Delta)^s)\varphi_2^\star,U_1-\varphi_1\rangle\,dt\\
        &=-\int_0^T\langle(-(-\Delta)^s\partial_t +(-\Delta)^s)U_2^\star,\varphi_1\rangle\,dt+\int_0^T\langle ((-\Delta)^s\partial_t +(-\Delta)^s)U_1,\varphi_2^\star\rangle\,dt\\
        &=-\int_0^T\langle[((-\Delta)^s\partial_t +(-\Delta)^s)U_2]^\star,\varphi_1\rangle\,dt+\int_0^T\langle ((-\Delta)^s\partial_t +(-\Delta)^s)U_1,\varphi_2^\star\rangle\,dt\\
        &=-\int_0^T\langle((-\Delta)^s\partial_t +(-\Delta)^s)U_2,\varphi_1^\star\rangle\,dt+\int_0^T\langle ((-\Delta)^s\partial_t +(-\Delta)^s)U_1,\varphi_2^\star\rangle\,dt.
        \end{split}
        \end{equation} 
         \normalsize
        By recalling the definition of the DN map associated with the viscous nonlocal wave equations \eqref{eq: for Q1} and \eqref{eq: for Q2star}, we get
        \begin{equation}
        \label{eq: caclc int identity diff pot}
             \begin{split}
        &\int_{\Omega_T}(Q_1-Q_2)(U_1-\varphi_1)(U_2-\varphi_2)^\star\,dxdt = \langle\Lambda_{Q_1}\varphi_1,\varphi_2^{\star}\rangle-\langle \Lambda_{Q_2^{\star}}\varphi_2,\varphi_1^{\star}\rangle.
        \end{split}
        \end{equation}
        On the one hand, choosing
        \[
            Q_1=Q_2=q_j
        \]
        in \eqref{eq: caclc int identity diff pot} ensures that
        \begin{equation}
        \label{eq: self adjointness DN map}
            \langle\Lambda_{q_j}\varphi_1,\varphi_2^\star\rangle=\langle \Lambda_{q_j^\ast}\varphi_2,\varphi_1^\star\rangle.
        \end{equation}
        On the other hand, taking 
        \[
            Q_1=q_1\text{ and }Q_2=q_2^\star
        \]
        in \eqref{eq: caclc int identity diff pot} yields
        \[
        \begin{split}
            &\int_{\Omega_T}(q_1-q_2^\star)(u_1-\varphi_1)(u_2-\varphi_2)^\star\,dxdt=\langle \Lambda_{q_1}\varphi_1,\varphi_2^\star\rangle-\langle\Lambda_{q_2}\varphi_2,\varphi_1^\star\rangle,
        \end{split}
        \]
        where we used that $U_j = u_j$, with $u_j$ denoting the solution to \eqref{eq: viscous wave for int id} for $j = 1,2$.
        Thus, using the identity \eqref{eq: self adjointness DN map}, we deduce that
        \[
        \begin{split}
            &\int_{\Omega_T}(q_1-q_2^\star)(u_1-\varphi_1)(u_2-\varphi_2)^\star\,dxdt =\langle(\Lambda_{q_1}-\Lambda_{q_2})\varphi_1,\varphi_2^\star\rangle.
        \end{split}
        \]
        Hence, we can conclude the proof.
    \end{proof}

    \begin{proof}[Proof of Theorem~\ref{Thm: main linear} for time-reversal invariant potentials] Throughout the proof we assume that $q_1, q_2$ are time-reversal invariant. Then, by Lemma~\ref{lemma: integral identity} and the condition \eqref{eq: cond DN map linear}, we get
        \[
            \int_{\Omega_T}(q_1-q_2)(u_1-\varphi_1)(u_2-\varphi_2)^{\star}dxdt=0,
        \]
        where $u_j$ is the unique solution of 
        \begin{equation}
            \begin{cases}
			(\partial_t^2 +(-\Delta)^s\partial_t+(-\Delta)^s + q_j)u = 0&\text{ in }\Omega_T \\
			u=\varphi_j  &\text{ in }(\Omega_e)_T,\\
			u(0)=0, \quad \partial_t u(0)=0 &\text{ in }\Omega.
	\end{cases}
        \end{equation}
        By the Runge approximation (Proposition~\ref{prop: runge}), for all $\Phi_1,\Phi_2\in C_c^{\infty}(\Omega_T)$ there exist sequences $v_k^{(1)}\in \mathcal{R}_{W_1}^{(1)}$ and $v_k^{(2)}\in \mathcal{R}_{W_2}^{(2)}$ such that
        \[
            v_k^{(1)}\to \Phi_1\text{ and }v_k^{(2)}\to \Phi_2\text{ in }L^2(0,T;\widetilde{H}^s(\Omega))
        \]
        as $k\to\infty$. Above we used the notation
        \begin{equation}
        \mathcal{R}_{W_j}^{(j)}=\{ u^{(j)}-\varphi\, : \, \varphi\in C^\infty_c ((W_j)_T) \},
		\end{equation}
		where $u^{(j)}\in \widetilde{W}_{ext}(0,T;\widetilde{H}^s(\Omega))$ is the unique solution to 
		\[
		\begin{cases}
			(\partial_t^2 +(-\Delta)^s\partial_t+(-\Delta)^s + q_j)u = 0 &\text{ in }\Omega_T \\
			u=\varphi  &\text{ in }(\Omega_e)_T,\\
			u(0)=0, \quad \partial_t u(0)=0 &\text{ in }\Omega.
		\end{cases}
		\] 
        Hence, we have
        \[
            \int_{\Omega_T}(q_1-q_2)v_k^{(1)}(v_k^{(2)})^{\star}dxdt=0
        \]
        for all $k\in\N$.
        By \eqref{eq: estimate a1}, we get
        \[
        \begin{split}
             &\left|\int_{\Omega_T}(q_1-q_2)v_k^{(1)}(v_k^{(2)})^{\star}dxdt-\int_{\Omega_T}(q_1-q_2)\Phi_1\Phi_2^{\star}dxdt\right|\to 0 
        \end{split}
        \]
        as $k\to\infty$. This shows that
        \[
            \int_{\Omega_T}(q_1-q_2)\Phi_1\Phi_2^{\star}dxdt=0.
        \]
        Hence, we may conclude that there holds $q_1=q_2$ a.e. in $\Omega_T$, which completes the proof.
    \end{proof}

    \begin{remark}
        Note that if one knows a priori that the potentials $q_j$ are bounded, then a Runge approximation in $L^2(\Omega_T)$ is enough to conclude the above uniqueness proof, but for lower regular potentials one needs the Runge approximation in $L^2(0,T;\widetilde{H}^s(\Omega))$.
    \end{remark}

     \begin{proof}[Proof of Theorem~\ref{Thm: main linear} for generic potentials]
        First, let us observe that \eqref{eq: cond DN map linear} ensures that the function 
\[
    w_{1,2} \vcentcolon= \partial_t u_{1,2} + u_{1,2} \in L^2(0,T;H^s(\mathbb{R}^n)),
\]
where \( u_{1,2} \vcentcolon= u_1 - u_2 \in H^1(0,T;H^s(\mathbb{R}^n)) \) and \(u_j\) denotes the solution to
\begin{equation}
\label{eq: viscous PDE 2nd proof}
\begin{cases}
    (\partial_t^2 + (-\Delta)^s \partial_t + (-\Delta)^s + q_j)u = 0 & \text{in } \Omega_T, \\
    u = \varphi & \text{in } (\Omega_e)_T, \\
    u(0)=0,\quad \partial_t u(0)=0 & \text{in } \Omega,
\end{cases}
\end{equation}
for a fixed exterior value \( \varphi \in C_c^\infty((W_1)_T) \), satisfies
\[
    (-\Delta)^s w_{1,2} = w_{1,2} \quad \text{in } (W_2)_T.
\]
Hence, the UCP for the fractional Laplacian \cite[Theorem~1.2]{GSU20} implies that \( w_{1,2} = 0 \) in \( \mathbb{R}^n_T \). By the same argument as in the proof of Proposition~\ref{prop: runge}, we then obtain \( u_{1,2} = 0 \) in \( \mathbb{R}^n_T \). Since \(u_j\) solves \eqref{eq: viscous PDE 2nd proof} and \(u_1 = u_2\) in \( \mathbb{R}^n_T \), we infer that
\[
    (q_1 - q_2)(u-\varphi) = 0 \quad \text{in } \Omega_T,
\]
where we set \( u \vcentcolon= u_1 = u_2 \). Repeating the arguments in the proof of Theorem~\ref{Thm: main linear} for time-reversal invariant potentials, we conclude that \( q_1 = q_2 \) in \( \Omega_T \).

     \end{proof}

    \section{Inverse problem for the nonlinear viscous wave equation}
    \label{IP nonlinear}

    In this section we study the inverse problem for the viscous wave equation with nonlinear perturbations. In Section~\ref{integral id nonlinear} we first introduce rigorously the DN map and then prove a suitable integral identity (Lemma~\ref{lem:integral identity}). Then finally in Section~\ref{unique determination of nonlinear perturb} we give the proof of Theorem~\ref{Thm: main nonlinear}.
	
	\subsection{An integral identity for the nonlinear problem}
    \label{integral id nonlinear}
	
     \begin{definition}[The DN map] 
    \label{DN map nonlinear case}
		Let $\Omega\subset\R^n$ be a bounded Lipschitz domain, $T>0$ and $s>0$ a non-integer. Suppose that we have given a nonlinearity $f$ satisfying Assumption~\ref{main assumptions on nonlinearities} and $f$ is $r+1$ homogeneous. Let $U_0\subset\widetilde{W}^s_{rest}((\Omega_e)_T),U_1\subset \widetilde{W}_{ext}(0,T;\widetilde{H}^s(\Omega))$ be the neighborhoods of Theorem~\ref{prop: differentiability of solution map} such that for any $\varphi\in U_0$ the problem 
        \begin{equation}
        \label{PDE def nonlinear DN map}
		\begin{cases}
			(\partial_t^2 +(-\Delta)^s\partial_t +(-\Delta)^s)u+f(u)= 0&\text{ in }\Omega_T\\
			u=\varphi&\text{ in }(\Omega_e)_T,\\
			u(0)=0,\quad \partial_t u(0)=0  &\text{ in }\Omega
		\end{cases}
		\end{equation}
        has a unique solution $u\in U_1$. Then we define the DN map $\Lambda_{f}$ related to
		\eqref{PDE def nonlinear DN map}
		by
		\begin{equation}
			\label{eq: modified DN map}
			\begin{split}
				\left\langle \Lambda_{f}\varphi_1,\varphi_2 \right\rangle \vcentcolon =	\int_{\R^n_T}(-\Delta)^{s/2}\,u(-\Delta)^{s/2}\varphi_2 \,dxdt+\int_{\R^n_T}(-\Delta)^{s/2}\partial_t u(-\Delta)^{s/2}\varphi_2 \,dxdt,
			\end{split}
		\end{equation}
		for all $\varphi\in U_0,\psi\in C_c^{\infty}((\Omega_e)_T)$, where $u\in U_1$ is the unique solution of \eqref{PDE def nonlinear DN map} with exterior condition $\varphi=\varphi_1$
		(see Theorem~\ref{prop: differentiability of solution map}). 
	\end{definition}

	\begin{lemma}\label{lem:integral identity}
		Let $\Omega\subset\R^n$ be a bounded Lipschitz domain, $T>0$ and $s>0$ a non-integer. Suppose that for $j=1,2$ we have given nonlinearities $f_j$ satisfying Assumption~\ref{main assumptions on nonlinearities} with $r>0$ and $a\in L^{\infty}(\Omega)$ and $f_j(0)=0$. Let $U^j_0\subset\widetilde{W}^s_{rest}((\Omega_e)_T),U^j_1\subset \widetilde{W}_{ext}(0,T;\widetilde{H}^s(\Omega))$ be the neighborhoods of Theorem~\ref{prop: differentiability of solution map} such that for any $\varphi\in U^j_0$ the problem 
        \begin{equation}
        \label{PDE integral identity nonlinear}
		\begin{cases}
			(\partial_t^2 +(-\Delta)^s\partial_t +(-\Delta)^s)u+f_j(u)= 0&\text{ in }\Omega_T\\
			u=\varphi&\text{ in }(\Omega_e)_T,\\
			u(0)=0,\quad \partial_t u(0)=0  &\text{ in }\Omega
		\end{cases}
		\end{equation}
        has a unique solution $u\in U_1^j$. Then for all exterior conditions $\varphi_1\in U_0^1\cap U_0^2\cap C_c^{\infty}((\Omega_e)_T)$ and $\varphi_2\in C_c^{\infty}((\Omega_e)_T)$ one has
		\begin{equation}\label{eq:integral identity}
			\left\langle \LC \Lambda_{f_1} - \Lambda_{f_2}\RC \varphi_1,\varphi_2^*\right\rangle \\
			=\int_{\Omega_T} (f_1(u_1^{(1)}) - f_2(u_1^{(2)}))(u_2-\varphi_2)^\star
			dxdt,
		\end{equation}
		where $u_1^{(j)}$ is the unique solution of \eqref{PDE integral identity nonlinear} with  $\varphi=\varphi_1$ and $u_2$ is the unique solution of the linear equation
		\begin{equation}\label{eq:linear eq for u_2}
			\begin{cases}
				\LC \p_t^2+(-\Delta)^s\partial_t+(-\Delta)^s\RC u=0, &\text{in } \Omega_T,\\
				u = \varphi_2, &\text{in } (\Omega_e)_T,\\
				u(0)=\p_t u(0) = 0,&\text{in } \Omega.
			\end{cases}    
		\end{equation}
	\end{lemma}
	\begin{proof}
		First of all note that $u=u_1^{(1)}-u_1^{(2)}\in \widetilde{W}(0,T;\widetilde{H}^s(\Omega))$ solves
        \begin{equation}
        \label{PDE difference of sols}
		\begin{cases}
			(\partial_t^2 +(-\Delta)^s\partial_t +(-\Delta)^s)u= -(f_1(u_1^{(1)})-f_2(u_1^{(2)}))&\text{ in }\Omega_T\\
			u=0&\text{ in }(\Omega_e)_T,\\
			u(0)=0,\quad \partial_t u(0)=0  &\text{ in }\Omega.
		\end{cases}
		\end{equation}
        Then the definition of the DN map \eqref{eq: modified DN map} implies
        \begin{equation}
        \label{eq: first id for integral identity}
        \begin{split}
            &\langle(\Lambda_{f_1}-\Lambda_{f_2})\varphi_1,\varphi_2^\star\rangle\\
            &=-\int_{0}^T\langle (-\Delta)^s(u_1^{(1)}-u_1^{(2)})+(-\Delta)^s\partial_t(u_1^{(1)}-u_1^{(2)}) ,(u_2-\varphi_2)^\star\rangle\,dt\\
            &\quad+\int_{0}^T\langle (-\Delta)^{s/2}(u_1^{(1)}-u_1^{(2)})+(-\Delta)^{s/2}\partial_t(u_1^{(1)}-u_1^{(2)}) ,(-\Delta)^{s/2}u_2^\star\rangle\,dt\\
            &=I_1+I_2.
        \end{split}
        \end{equation}
        As $u_1^{(1)}-u_1^{(2)}$ solves \eqref{PDE difference of sols}, the identity \eqref{eq: integration by parts 2nd order der} together with the fact that $u_2$ is a solution of \eqref{eq:linear eq for u_2} implies
        \begin{equation}
        \begin{split}
            I_1&=\int_0^T\langle \partial_t^2(u_1^{(1)}-u_1^{(2)})+(f_1(u_1^{(1)})-f_2(u_1^{(2)})),(u_2-\varphi_2)^\star\rangle\,dt\\
            &=\int_0^T(\langle \partial_t^2(u_2-\varphi_2)^\star,(u_1^{(1)}-u_1^{(2)})\rangle+\langle(f_1(u_1^{(1)})-f_2(u_1^{(2)})),(u_2-\varphi_2)^\star\rangle)\,dt\\
            &=-\int_0^T\langle (-(-\Delta)^s\partial_t +(-\Delta)^s )(u_2-\varphi_2)^\star,u_1^{(1)}-u_1^{(2)}\rangle\,dt \\
            &\quad -\int_0^T\langle (\partial_t^2 -(-\Delta)^s\partial_t +(-\Delta)^s)\varphi_2^\star,u_1^{(1)}-u_1^{(2)}\rangle \,dt\\
            &\quad +\int_0^T \langle(f_1(u_1^{(1)})-f_2(u_1^{(2)})),(u_2-\varphi_2)^\star\rangle)\,dt
        \end{split}
        \end{equation}
        Taking into account the support condition of $\varphi_2$, we deduce that
        \begin{equation}
        \label{eq: I1}
        \begin{split}
            I_1&=-\int_0^T\langle (-(-\Delta)^s\partial_t +(-\Delta)^s )u_2^\star,u_1^{(1)}-u_1^{(2)}\rangle\,dt \\
            &\quad +\int_0^T \langle(f_1(u_1^{(1)})-f_2(u_1^{(2)})),(u_2-\varphi_2)^\star\rangle)\,dt.
        \end{split}
        \end{equation}
        On the other hand, using integration by parts, the integral $I_2$ is given by
        \begin{equation}
            \label{eq: I2}
            \begin{split}
                I_2&=\int_{0}^T\langle (-(-\Delta)^{s/2}\partial_t+(-\Delta)^{s/2})u_2^\star ,(-\Delta)^{s/2}(u_1^{(1)}-u_1^{(2)})\rangle\,dt
            \end{split}
        \end{equation}
        Summing up \eqref{eq: I1} and \eqref{eq: I2}, we deduce from \eqref{eq: first id for integral identity} the desired identity \eqref{eq:integral identity} and can conclude the proof.
	\end{proof}
	
	\subsection{Unique determination of nonlinear perturbations}
    \label{unique determination of nonlinear perturb}
	
	    \begin{proof}[Proof of Theorem~\ref{Thm: main nonlinear}]
        Let $\boldsymbol{\varepsilon}=(\varepsilon_0,\varepsilon_1)$ and define
        \begin{equation}
        \label{eq: test functions uniqueness proof}
            \varphi_1^{\Vareps}=\eps_0\psi_0+\eps_1\psi_1
        \end{equation}
        for some $\psi_j\in C_c^{\infty}((W_1)_T)$. If $|\Vareps|\ll 1$, then by the integral identity \eqref{eq:integral identity} with $\varphi_1=\varphi_1^{\Vareps}$ of Lemma~\ref{lem:integral identity} we have
        \begin{equation}
            \label{eq:integral identity uniqueness proof}
			\left\langle \LC \Lambda_{f_1} - \Lambda_{f_2}\RC \varphi^{\Vareps}_1,\varphi_2^*\right\rangle \\
			=\int_{\Omega_T} (f_1(u_{\Vareps}^{(1)}) - f_2(u_{\Vareps}^{(2)}))(u_2-\varphi_2)^\star
			dxdt
        \end{equation}
        for $\varphi_2\in C_c^{\infty}((W_2)_T)$, where $u_{\Vareps}^{(j)}$ solves 
        \begin{equation}
        \label{eq: approx sols uniqueness proof}
        \begin{cases}
            (\partial_t^2 +(-\Delta)^s\partial_t +(-\Delta)^s)u+f_j(u)= 0&\text{ in }\Omega_T\\
			u=\varphi_1^{\Vareps}&\text{ in }(\Omega_e)_T,\\
			u(0)=0,\quad \partial_t u(0)=0  &\text{ in }\Omega
        \end{cases}
        \end{equation}
        and $u_2$ solves
        \begin{equation}
            \label{eq:linear eq for u_2 uniqueness proof}
			\begin{cases}
				\LC \p_t^2+(-\Delta)^s\partial_t+(-\Delta)^s\RC u=0, &\text{in } \Omega_T,\\
				u = \varphi_2, &\text{in } (\Omega_e)_T,\\
				u(0)=\p_t u(0) = 0,&\text{in } \Omega.
			\end{cases}  
        \end{equation}
    By Theorem~\ref{prop: differentiability of solution map}, we know that the solution map $U_0^j\ni \varphi\mapsto u=S_j(\varphi)$ associated to
    \begin{equation}
        \label{eq: diff sol map uniqueness proof}
        \begin{cases}
            (\partial_t^2 +(-\Delta)^s\partial_t +(-\Delta)^s)u+f_j(u)= 0&\text{ in }\Omega_T\\
			u=\varphi&\text{ in }(\Omega_e)_T,\\
			u(0)=0,\quad \partial_t u(0)=0  &\text{ in }\Omega
        \end{cases}
        \end{equation}
    is $C^1$ as a map from $U_0^j$ to $\widetilde{W}_{ext}(0,T;\widetilde{H}^s(\Omega))$. In particular, we see that
    \begin{equation}
    \label{eq: derivative of solution map}
        v_j=\left.\partial_{\epsilon}\right|_{\epsilon=0}S_j(\rho+\epsilon\eta)\in \widetilde{W}_{ext}(0,T;\widetilde{H}^s(\Omega))
    \end{equation}
    exists for any $\rho\in U_0^j$ and $\eta\in C_c^{\infty}((\Omega_e)_T)$. Moreover, from the proof of Theorem~\ref{prop: differentiability of solution map} it follows that $v_j$ solves
    \begin{equation}
    \label{eq: equation for derivative}
         \begin{cases}
            (\partial_t^2 +(-\Delta)^s\partial_t +(-\Delta)^s+\partial_{\tau}f_j(S_j(\rho)))v= 0&\text{ in }\Omega_T\\
			v=\eta&\text{ in }(\Omega_e)_T,\\
			v(0)=0,\quad \partial_t v(0)=0  &\text{ in }\Omega.
        \end{cases}
    \end{equation}
    Using these observations we deduce from \eqref{eq:integral identity uniqueness proof} that there holds
    \begin{equation}
    \label{eq: derivative in eps 0 of DN map}
        \begin{split}
            &\left.\partial_{\varepsilon_1}\right|_{\varepsilon_1=0}\langle (\Lambda_{f_1}-\Lambda_{f_2})\varphi_1^{\Vareps},\varphi_2^\star\rangle\\
            &=\langle \partial_{\tau}f_1(u_{\varepsilon_0}^{(1)})v_{\eps_0,\psi_1}^{(1)}-\partial_{\tau}f_2(u_{\eps_0}^{(2)})v_{\eps_0,\psi_1}^{(2)},(u_2-\varphi_2^{\star})\rangle,
        \end{split}
    \end{equation}
    where the right hand side is the duality pairing between $L^2(0,T;\widetilde{H}^s(\Omega))$ and $L^2(0,T;H^{-s}(\Omega))$ and we set
    \begin{equation}
    \label{eq: definition of functions}
    \begin{split}
         u_{\eps_0}^{(j)}&=\lim_{\eps_1\to 0}u_{\Vareps}^{(j)}\\
         v_{\eps_0,\psi_1}^{(j)}&=\left.\partial_{\eps_1}\right|_{\eps_1=0}u_{\Vareps}^{(j)}.
    \end{split}
    \end{equation}
    By the continuity of the solution map and \eqref{eq: equation for derivative}, it follows that $u_{\eps_0}^{(j)}$ and $v_{\eps_0,\psi_1}^{(j)}$ solve
    \begin{equation}
    \label{eq: PDE for u eps0}
        \begin{cases}
            (\partial_t^2 +(-\Delta)^s\partial_t +(-\Delta)^s)u+f_j(u)= 0&\text{ in }\Omega_T\\
			u=\eps_0\psi_0&\text{ in }(\Omega_e)_T,\\
			u(0)=0,\quad \partial_t u(0)=0  &\text{ in }\Omega
        \end{cases}
    \end{equation}
    and 
    \begin{equation}
     \label{eq: PDE for v eps0 psi1}
        \begin{cases}
            (\partial_t^2 +(-\Delta)^s\partial_t +(-\Delta)^s+\partial_{\tau}f_j(u_{\eps_0}^{(j)}))v= 0&\text{ in }\Omega_T\\
			v=\psi_1&\text{ in }(\Omega_e)_T,\\
			v(0)=0,\quad \partial_t v(0)=0  &\text{ in }\Omega,
        \end{cases}
    \end{equation}
    respectively. Next note that by the homogenity of $\partial_\tau f_j$ and arguing as above, one has
    \begin{equation}
    \label{eq: convergence of first term}
        \eps_0^{-1}u_{\eps_0}^{(j)}\to v_0=\left.\partial_{\eps_0}\right|_{\eps_0=0}u_{\eps_0}^{(j)}
    \end{equation}
    as $\eps_0\to 0$ in $\widetilde{W}_{ext}(0,T;\widetilde{H}^s(\Omega))$ and in particular in $L^q(0,T;H^s(\R^n))$ for any $1\leq q\leq \infty$. Moreover, $v_0$ is the unique solution of
    \begin{equation}
    \label{eq: PDE for v0}
        \begin{cases}
            (\partial_t^2 +(-\Delta)^s\partial_t +(-\Delta)^s)v= 0&\text{ in }\Omega_T\\
			v=\psi_0&\text{ in }(\Omega_e)_T,\\
			v(0)=0,\quad \partial_t v(0)=0  &\text{ in }\Omega.
        \end{cases}
    \end{equation}
    Next, we show the following assertion.
    \begin{claim}
    \label{claim: convergence of second term}
        Let $w_{\eps_0,\psi_1}^{(j)}$ be defined by
        \begin{equation}
        \label{eq: def of w function}
            w_{\eps_0,\psi_1}^{(j)}=v_{\eps_0,\psi_1}^{(j)}-\psi_1\in\widetilde{W}(0,T;\widetilde{H}^s(\Omega)).
        \end{equation}
        Then there exists $w_1\in \widetilde{W}(0,T;\widetilde{H}^s(\Omega))$ such that 
        \begin{enumerate}[(i)]
            \item\label{weak conv} $w_{\eps_0,\psi_1}^{(j)}\weak w_1$ in $\widetilde{W}(0,T;\widetilde{H}^s(\Omega))$,
            \item\label{strong conv} $w_{\eps_0,\psi_1}^{(j)}\to w_1$ in $C([0,T];\widetilde{H}^t(\Omega))$ for any $0\leq t<s$.
        \end{enumerate}
        Moreover, $v_1=w_1+\psi_1$ is the unique solution of
        \begin{equation}
            \label{eq: PDE for v1}
        \begin{cases}
            (\partial_t^2 +(-\Delta)^s\partial_t +(-\Delta)^s)v= 0&\text{ in }\Omega_T\\
			v=\psi_1&\text{ in }(\Omega_e)_T,\\
			v(0)=0,\quad \partial_t v(0)=0  &\text{ in }\Omega.
        \end{cases}
        \end{equation}
    \end{claim}
    \begin{proof}[Proof of Claim~\ref{claim: convergence of second term}]
        Let us start by observing that $w_{\eps_0,\psi_1}^{(j)}$ is the unique solution of
        \begin{equation}
            \label{eq: PDE for weps0}
        \begin{cases}
            (\partial_t^2 +(-\Delta)^s\partial_t +(-\Delta)^s+\partial_{\tau}f_j(u_{\eps_0}^{(j)}))w= -((-\Delta)^s\partial_t+(-\Delta)^s)\psi_1&\text{ in }\Omega_T\\
			w=0&\text{ in }(\Omega_e)_T,\\
			w(0)=0,\quad \partial_t w(0)=0  &\text{ in }\Omega.
        \end{cases}
        \end{equation}
        Now, we show that the function
        \[
            q: =\partial_{\tau}f_j(u_{\eps_0}^{(j)})\in L^1_{loc}(\Omega_T)
        \]
        satisfies the conditions \ref{integrability of q} and \ref{continuity of q} of Theorem~\ref{thm: well-posedness linear}. 

        \noindent \ref{integrability of q}: If $2s\geq n$, then the integrability is clear as $H^s(\R^n)\hookrightarrow L^p(\R^n)$ for any $2\leq p<\infty$. Hence, we can assume without loss of generality that $2s<n$. In this case the condition
        \[
            0<r\leq \frac{2s}{n-2s}
        \]
        guarantees that
        \[
            \frac{n}{s}\leq \frac{2n}{r(n-2s)}.
        \]
        Therefore, we may estimate
        \begin{equation}
        \label{eq: L n/s bound}
        \begin{split}
             \|\partial_\tau f_j(u_{\eps_0}^{(j)})\|_{L^{n/s}(\Omega)}&\lesssim \|\partial_\tau f_j(u_{\eps_0}^{(j)})\|_{L^{\frac{2n}{r(n-2s)}}(\Omega)}\lesssim \|u_{\eps_0}^{(j)}\|^r_{L^{\frac{2n}{n-2s}}(\Omega)}\\
             &\lesssim \|u_{\eps_0}^{(j)}\|^r_{H^s(\R^n)}<\infty.
        \end{split}
        \end{equation}
        Here, we are using that $\partial_\tau f_j (x,\tau)$ is $r$ homogeneous in the second variable, since $f_j$ is $(r+1)$ homogeneous.
        As $u_{\eps_0}^{(j)}\in L^{\infty}(0,T;H^s(\R^n))$ this shows that $\partial_\tau f_j(u_{\eps_0}^{(j)})\in L^{\infty}(0,T;L^{n/s}(\Omega))$ as we wanted to prove.\\

        \noindent \ref{continuity of q}: Next, note that
        \[
            u_{\eps_0}^{(j)}=u_{\eps_0}^{(j)}-\eps_0\psi_0\text{ in }\Omega_T.
        \]
        By Lemma~\ref{lemma: embeddings} and the Sobolev embedding, we know that 
        \begin{equation}
        \label{eq: properties of ueps0}
            u_{\eps_0}^{(j)}-\eps_0\psi_0\in C([0,T];\widetilde{H}^s(\Omega))\hookrightarrow C([0,T];L^{\bar{p}}(\Omega)).
        \end{equation}
        for all $\bar{p}$ satisfying
        \[
            \begin{cases}
                1\leq \bar{p}\leq \frac{2n}{n-2s},&\text{ for }2s<n\\
                1\leq \bar{p}<\infty,&\text{ for }2s=n\\
                1\leq \bar{p}\leq \infty,&\text{ for }2s>n.
            \end{cases}
        \]
        But then the continuity of
        \[
            t\mapsto \int_{\Omega}\partial_\tau f_j(u_{\eps_0}^{(j)})\varphi\,dx,
        \]
        for fixed $\varphi\in C_c^{\infty}(\Omega)$, is an immediate consequence of H\"older's inequality, the fact that $\partial_\tau f (u)$ is continuous as a map from $L^{r+2}(\Omega)$ to $L^{\frac{r+2}{r}}(\Omega)$ and that there holds
        \[
            2+r<\frac{2n}{n-2s}\text{ for }2s<n.
        \]
        This establishes the condition \ref{continuity of q}.\\

        Therefore, $q=\partial_{\tau}f_j(u_{\eps_0}^{(j)})$ satisfies all necessary conditions in Theorem~\ref{thm: well-posedness linear}, where
        \begin{equation}
        \label{eq: range of p uniqueness nonlinear}
            \begin{cases}
                p=n/s,&\text{ for }2s<n\\
                2<p<\infty,&\text{ for }2s=n\\
                2\leq p\leq \infty,&\text{ for }2s>n,
            \end{cases}
        \end{equation}
        and we can apply the energy identity of that theorem to obtain
        \small
        \begin{equation}
			\label{eq: energy identity w}
			\begin{split}
				& \|\partial_t w_{\eps_0,\psi_1}^{(j)}(t) \|_{L^2(\Omega)}^2+ \|(-\Delta)^{s/2}w_{\eps_0,\psi_1}^{(j)}(t)\|_{L^2(\R^n)}^2+ 2\|(-\Delta)^{s/2}\partial_t w_{\eps_0,\psi_1}^{(j)}\|_{L^2(\R^n_t)}^2\\
				&=-2\int_0^t\langle ((-\Delta)^s\partial_t+(-\Delta)^s)\psi_1(\sigma),\partial_t w_{\eps_0,\psi_1}^{(j)} (\sigma)\rangle\,d\sigma-2\langle \partial_{\tau}f_j(u_{\eps_0}^{(j)}) w_{\eps_0,\psi_1}^{(j)},\partial_t w_{\eps_0,\psi_1}^{(j)} \rangle_{L^2(\Omega_t)}.
			\end{split}
		\end{equation}
        \normalsize
        By \eqref{eq: estimate a1}, we know that there holds
        \begin{equation}
        \label{eq: uniform bound potential}
            \left|\langle \partial_{\tau}f_j(u_{\eps_0}^{(j)}) u,v\rangle_{L^2(\Omega)} \right|\leq C\|\partial_{\tau}f_j(u_{\eps_0}^{(j)}(t))\|_{L^p(\Omega)}\|u\|_{\widetilde{H}^s(\Omega)}\|v\|_{L^2(\Omega)}
        \end{equation}
        for all $u,v\in\widetilde{H}^s(\Omega)$. Hence, the last term in the second line of \eqref{eq: energy identity w} can be estimate as
        \begin{equation}
        \label{eq: estimate last term in uniqueness f}
            \begin{split}
                &\left|\langle \partial_{\tau}f_j(u_{\eps_0}^{(j)}) w_{\eps_0,\psi_1}^{(j)},\partial_t w_{\eps_0,\psi_1}^{(j)} \rangle_{L^2(\Omega_t)}\right|\\
                &\quad \leq C\int_0^t \|\partial_{\tau}f_j(u_{\eps_0}^{(j)}(\sigma))\|_{L^p(\Omega)}(\|w_{\eps_0,\psi_1}^{(j)}(\sigma)\|_{\widetilde{H}^s(\Omega)}^2+\|\partial_t w_{\eps_0,\psi_1}^{(j)}(\sigma)\|_{L^2(\Omega)}^2)\,d\sigma
            \end{split}
        \end{equation}
        for some constant $C>0$ independent of $T$. On the other hand, the first term in the second line of \eqref{eq: energy identity w} can be estimated as
        \begin{equation}
        \label{eq: estimate last term in uniqueness ext}
            \begin{split}
                &\left|\int_0^t\langle ((-\Delta)^s\partial_t+(-\Delta)^s)\psi_1(\sigma),\partial_t w_{\eps_0,\psi_1}^{(j)} (\sigma)\rangle\,d\sigma\right|\\
                &\leq C\int_0^t\|(-\Delta)^s\partial_t+(-\Delta)^s)\psi_1(\sigma)\|_{H^{-s}(\Omega)}^2+\int_0^t\|\partial_t w_{\eps_0,\psi_1}^{(j)}(\sigma)\|_{\widetilde{H}^s(\Omega)}^2\,d\sigma
            \end{split}
        \end{equation}
        Note that we can absorb the last term on the left hand side of \eqref{eq: energy identity w}. Then combining \eqref{eq: energy identity w},   \eqref{eq: estimate last term in uniqueness f} and \eqref{eq: estimate last term in uniqueness ext}, we see that the function
        \begin{equation}
        \label{eq: function for gronwall}
            \Phi(t)=\|\partial_t  w_{\eps_0,\psi_1}^{(j)}(t)\|_{L^2(\Omega)}^2+\| (-\Delta)^{s/2} w_{\eps_0,\psi_1}^{(j)}(t)\|_{L^2(\R^n)}+\|(-\Delta)^{s/2}\partial_t  w_{\eps_0,\psi_1}^{(j)}\|_{L^2(\R^n_t)}^2
        \end{equation}
        satisfies $\Phi\in L^{\infty}((0,T))$ and
        \small
        \begin{equation}
        \label{eq: estimate for gronwall}
            \Phi(t)\leq C\left(\int_0^t \|\partial_{\tau}f_j(u_{\eps_0}^{(j)}(\sigma))\|_{L^p(\Omega)}\Phi(\sigma)\,d\sigma+\int_0^t\|(-\Delta)^s\partial_t+(-\Delta)^s)\psi_1(\sigma)\|_{H^{-s}(\Omega)}^2 d\sigma\right).
        \end{equation}
        \normalsize
        Now, using $\partial_\tau f_j(u_{\eps_0}^{(j)})\in L^{\infty}(0,T;L^{p}(\Omega))$, we deduce from Gronwall's inequality the estimate
        \begin{equation}
        \label{eq: estimate from Gronwall}
            \Phi(t)\leq C\|(-\Delta)^s\partial_t+(-\Delta)^s)\psi_1\|_{L^2(0,T;H^{-s}(\Omega))}^2e^{C\|\partial_\tau f_j(u_{\eps_0}^{(j)})\|_{L^1(0,T;L^p(\Omega))}}
        \end{equation}
        for a.e. $0\leq t\leq T$.

        Next recall by \eqref{eq: cond f and partial f} that the map
        \begin{equation}
        \label{eq: continuity of der of f}
            \partial_{\tau}f\colon L^q(0,T;H^s(\R^n))\to L^{q/r}(0,T;L^{n/s}(\Omega))
        \end{equation}
        is continuous for any $r\leq q<\infty$. As the solution map is continuous and $f_j(0)=0$, we deduce that
        \[
            u_{\eps_0}^{(j)}\to 0\text{ in }L^{\infty}(0,T;H^s(\R^n))
        \]
        as $\eps_0\to 0$. This ensures that
        \begin{equation}
        \label{eq: conv of potential uniqueness}
            \partial_{\tau}f_j(u_{\eps_0}^{(j)})\to 0\text{ in }L^{q/r}(0,T;L^{n/s}(\Omega))\hookrightarrow L^1(0,T;L^{n/s}(\Omega))
        \end{equation}
        as $\eps_0\to 0$, for any $r\leq q<\infty$. Hence, by  \eqref{eq: estimate from Gronwall} we achieve that
        \begin{equation}
        \label{eq: uniform bound}
            w_{\eps_0,\psi_1}^{(j)}\text{ is uniformly bounded in }\widetilde{W}(0,T;\widetilde{H}^s(\Omega)).
        \end{equation}
        Note that the uniform bound of $\partial_t^2 w_{\eps_0,\psi_1}^{(j)}$ in $L^2(0,T;H^{-s}(\Omega))$ comes from the PDE \eqref{eq: PDE for weps0}, the uniform bound of $w_{\eps_0,\psi_1}$ in $H^1(0,T;\widetilde{H}^s(\Omega))$ and the estimate \eqref{eq: L n/s bound} (with similar estimate in the range $2s\geq n$). By the usual embeddings we get
        \begin{enumerate}[(a)]
            \item\label{boundedness 1} $w_{\eps_0,\psi_1}^{(j)}$ is uniformly bounded in $L^{\infty}(0,T;\widetilde{H}^s(\Omega))$,
            \item\label{boundedness 2} $\partial_t w_{\eps_0,\psi_1}^{(j)}$ is uniformly bounded in $C([0,T];L^2(\Omega))$.
        \end{enumerate}
        Next, recall that $\widetilde{H}^s(\Omega)\hookrightarrow \widetilde{H}^t(\Omega)\hookrightarrow L^2(\Omega)$ for any $0<t<s$, where the first embedding is compact. Using \ref{boundedness 1}, \ref{boundedness 2} and the Aubin--Lions lemma (\cite[Corollary~4]{Simon}) we see that 
        \begin{equation}
        \label{eq: rel compactness}
            w_{\eps_0,\psi_1}^{(1)}\text{ is relatively compact in }C([0,T];\widetilde{H}^t(\Omega))
        \end{equation}
        for any $0<t<s$. Thus, we deduce from \eqref{eq: uniform bound} and \eqref{eq: rel compactness} that there exists $w_1^{(j)}\in\widetilde{W}(0,T;\widetilde{H}^s(\Omega))$ and a subsequence of $w_{\eps_0,\psi_1}^{(j)}$, $\eps_0>0$, such that
        \begin{enumerate}[(I)]
            \item\label{conv 1} $w_{\eps_0,\psi_1}^{(j)}\weak w_1^{(j)}$ in $\widetilde{W}(0,T;\widetilde{H}^s(\Omega))$ as $\eps_0\to 0$, 
            \item\label{conv 2} $w_{\eps_0,\psi_1}^{(j)}\to w_1^{(j)}$ in $C([0,T];\widetilde{H}^t(\Omega))$ for all $0<t<s$ as $\eps_0\to 0$.
        \end{enumerate}
        The convergence \ref{conv 1} clearly implies
        \begin{equation}
        \label{eq: convergence of 2nd and 3rd term in PDE}
        \begin{split}
            (-\Delta)^sw_{\eps_0,\psi_1}^{(j)}&\weak (-\Delta)^s w_1^{(j)}\text{ in }L^2(0,T;H^{-s}(\Omega))\\
            (-\Delta)^s\partial_t w_{\eps_0,\psi_1}^{(j)}&\weak (-\Delta)^s \partial_t w_1^{(j)}\text{ in }L^2(0,T;H^{-s}(\Omega))
        \end{split}
        \end{equation}
        as $\eps_0\to 0$. On the other hand, by choosing $q$ sufficiently large in \eqref{eq: conv of potential uniqueness} we see that
        \begin{equation}
        \label{eq: continuity estimate uniqueness proof for later}
        \begin{split}
           & \left|\int_{\Omega_T} \partial_\tau f_j(u_{\eps_0}^{(j)})w_{\eps_0,\psi_1}^{(j)} v\,dxdt\right|\leq \int_{0}^T\|\partial_\tau f_j(u_{\eps_0}^{(j)})\|_{L^{n/s}(\Omega)}\|v\|_{\widetilde{H}^s(\Omega)}\|w_{\eps_0,\psi_1}^{(j)}\|_{L^2(\Omega)}\,dt\\
            &\quad \leq \|w_{\eps_0,\psi_1}^{(j)}\|_{L^{\infty}(0,T;L^2(\Omega))}\|\partial_\tau f_j(u_{\eps_0}^{(j)})\|_{L^2(0,T;L^{n/s}(\Omega))}\|v\|_{L^2(0,T;\widetilde{H}^s(\Omega))}
        \end{split}
        \end{equation}
        for all $v\in L^2(0,T;\widetilde{H}^s(\Omega))$. As the first factor is uniformly bounded, we get the convergence
        \begin{equation}
        \label{eq: convergence nonlinearity}
            \partial_\tau f_j(u_{\eps_0}^{(j)})w_{\eps_0,\psi_1}^{(j)}\weak 0\text{ in }L^2(0,T;H^{-s}(\Omega))
        \end{equation}
        as $\eps_0\to 0$. Again a similar argument can be used in the cases $2s\geq n$ to obtain this convergence.
        Using \eqref{eq: convergence of 2nd and 3rd term in PDE} and \eqref{eq: convergence nonlinearity}, we can pass to the limit in the weak formulation of \eqref{eq: PDE for weps0} and see that $w_1^{(j)}$ solves
        \begin{equation}
        \label{eq: PDE for w1}
            (\partial_t^2 +(-\Delta)^s\partial_t +(-\Delta)^s)w= -((-\Delta)^s\partial_t+(-\Delta)^s)\psi_1
        \end{equation}
        in $\Omega_T$. Additionally, from the trace theorem, we infer that
        \begin{equation}
        \label{eq: initial cond}
            w_1^{(j)}(0)=\partial_t w_1^{(j)}(0)=0.
        \end{equation}
        Thus, $w_1^{(j)}$ is the unique solution of 
        \[
        \begin{cases}
            (\partial_t^2 +(-\Delta)^s\partial_t +(-\Delta)^s)w= -((-\Delta)^s\partial_t+(-\Delta)^s)\psi_1&\text{ in }\Omega_T\\
			w=0&\text{ in }(\Omega_e)_T,\\
			w(0)=0,\quad \partial_t w(0)=0  &\text{ in }\Omega.
        \end{cases}
        \]
        and in particular is independent of $j$. As the above analysis works for any subsequence of $w_{\eps_0,\psi_1}^{(j)}$, one also sees that the whole sequence needs to converge to $w_1$.
        Therefore, we may conclude from \eqref{eq: PDE for w1} and \eqref{eq: initial cond} that $v$ is the unique solution of  \eqref{eq: PDE for v1} and this finishes the proof of the Claim~\ref{claim: convergence of second term}.
    \end{proof}
    Next, let us recall that by \eqref{eq: derivative in eps 0 of DN map} and  \eqref{eq: cond DN map nonlinear}, we have
    \begin{equation}
         0=\langle \partial_{\tau}f_1(u_{\varepsilon_0}^{(1)})v_{\eps_0,\psi_1}^{(1)}-\partial_{\tau}f_2(u_{\eps_0}^{(2)})v_{\eps_0,\psi_1}^{(2)},(u_2-\varphi_2)^{\star}\rangle.
    \end{equation}
    Multiplying this identity by $\epsilon_0^{-r}$ and using the $r$ homogeneity of $\partial_\tau f_j(u)$, we get
    \begin{equation}
    \label{eq: identity for uniqueness}
        0=\langle \partial_{\tau}f_1(\eps_0^{-1}u_{\varepsilon_0}^{(1)})v_{\eps_0,\psi_1}^{(1)}-\partial_{\tau}f_2(\eps_0^{-1}u_{\eps_0}^{(2)})v_{\eps_0,\psi_1}^{(2)},(u_2-\varphi_2)^{\star}\rangle.
    \end{equation}
    By \eqref{eq: convergence of first term} and \eqref{eq: continuity of der of f}, we know that
    \begin{equation}
    \label{eq: conv nonlinear term}
        \partial_{\tau}f_j(\eps_0^{-1}u_{\varepsilon_0}^{(j)})\to \partial_\tau f_j(v_0)\text{ in }L^{q/r}(0,T;L^{n/s}(\Omega))
    \end{equation}
    as $\eps_0\to 0$ for any $q\geq \max (1,r)$. If we choose $q$ such that $q\geq \max(1,2r)$, then the computation in \eqref{eq: continuity estimate uniqueness proof for later} shows that
    \[
        \begin{split}
           & \left|\int_{\Omega_T} \partial_\tau f_j(\eps_0^{-1}u_{\eps_0}^{(j)})w_{\eps_0,\psi_1}^{(j)} \eta\,dxdt\right|\\
           &\leq \|w_{\eps_0,\psi_1}^{(j)}\|_{L^{\infty}(0,T;L^2(\Omega))}\|\partial_\tau f_j(\eps_0^{-1}u_{\eps_0}^{(j)})\|_{L^2(0,T;L^{n/s}(\Omega))}\|\eta\|_{L^2(0,T;\widetilde{H}^s(\Omega))}\\
           &\leq C\|w_{\eps_0,\psi_1}^{(j)}\|_{L^{\infty}(0,T;L^2(\Omega))}\|\partial_\tau f_j(\eps_0^{-1}u_{\eps_0}^{(j)})\|_{L^{q/r}(0,T;L^{n/s}(\Omega))}\|\eta\|_{L^2(0,T;\widetilde{H}^s(\Omega))}
        \end{split}
    \]
    for any $\eta\in L^2(0,T;\widetilde{H}^s(\Omega))$. Using this estimate, the convergence \ref{strong conv} and \eqref{eq: conv nonlinear term}, we deduce that
    \begin{equation}
    \label{eq: final conv}
        \int_{\Omega_T} \partial_\tau f_j(\eps_0^{-1}u_{\eps_0}^{(j)})w_{\eps_0,\psi_1}^{(j)} \eta\,dxdt\to \int_{\Omega_T} \partial_\tau f_j(v_0)w_1 \eta\,dxdt
    \end{equation}
    for any $\eta\in L^2(0,T;\widetilde{H}^s(\Omega))$. In particular, \eqref{eq: final conv} and the splitting $v_{\eps_0,\psi_1}^{(j)}=w_{\eps_0,\psi_1}^{(j)}+\psi_1$ allows us to pass to the limit in \eqref{eq: identity for uniqueness}, which gives
    \begin{equation}
    \label{eq: identity for limit in uniqueness proof}
         0=\langle (\partial_{\tau}f_1(v_0)-\partial_{\tau}f_2(v_0))v_1,(u_2-\varphi_2)^{\star}\rangle.
    \end{equation}
    Now, let $\Psi_j\in C_c^{\infty}(\Omega_T)$, $j=0,1,2$, be given functions and choose according to the Runge approximation (Proposition~\ref{prop: runge}) the following sequences:
    \begin{enumerate}[(A)]
         \item $v^k_0-\psi_0^k\in L^2(0,T;\widetilde{H}^s(\Omega))$, $k\in\N$, where $\psi_0^k\in C_c^{\infty}((W_1)_T)$, $v^k_0$ is the unique solution of 
        \[
            \begin{cases}
            (\partial_t^2 +(-\Delta)^s\partial_t +(-\Delta)^s)v= 0&\text{ in }\Omega_T\\
			v=\psi^k_0&\text{ in }(\Omega_e)_T,\\
			v(0)=0,\quad \partial_t v(0)=0  &\text{ in }\Omega
        \end{cases}
        \]
        and $v^k_0-\psi_0^k\to \Psi_0$ in $L^2(0,T;\widetilde{H}^s(\Omega))$.
        \item $v^k_1-\psi_1^k\in L^2(0,T;\widetilde{H}^s(\Omega))$, $k\in\N$, where $\psi_1^k \in C_c^{\infty}((W_1)_T)$, $v^k_1$ is the unique solution of 
        \[
            \begin{cases}
            (\partial_t^2 +(-\Delta)^s\partial_t +(-\Delta)^s)v= 0&\text{ in }\Omega_T\\
			v=\psi^k_1&\text{ in }(\Omega_e)_T,\\
			v(0)=0,\quad \partial_t v(0)=0  &\text{ in }\Omega
        \end{cases}
        \]
        and $v^k_1-\psi_1^k\to \Psi_1$ in $L^2(0,T;\widetilde{H}^s(\Omega))$ as well as $v^k_1\to \Psi_1$ in $L^2(\Omega_T)$.
        \item $u^k_2-\varphi_2^k\in L^2(0,T;\widetilde{H}^s(\Omega))$, $k\in\N$, where $\varphi_2^k\in C_c^{\infty}((W_2)_T)$, $u^k_2$ is the unique solution of 
        \[
        \begin{cases}
				\LC \p_t^2+(-\Delta)^s\partial_t+(-\Delta)^s\RC u=0, &\text{in } \Omega_T,\\
				u = \varphi_2^k, &\text{in } (\Omega_e)_T,\\
				u(0)=\p_t u(0) = 0,&\text{in } \Omega
			\end{cases}  
        \]
        and $u^k_2-\varphi_2^k\to \Psi_2$ in $L^2(0,T;\widetilde{H}^s(\Omega))$.
    \end{enumerate}
    As seen in the beginning of the proof of Claim~\ref{claim: convergence of second term}, we know that if $v\in L^{\infty}(0,T;H^s(\R^n))$, then 
    \[
        q_j=\partial_{\tau}f_j(v)\in L^{\infty}(0,T;L^p(\Omega))
    \]
    for some $p$ satisfying the restrictions given in equation \eqref{eq: range of p uniqueness nonlinear}. Moreover, from \eqref{eq: estimate a1} we deduce that there holds
    \begin{equation}
    \label{eq: estimate for difference}
        \begin{split}
            \left|\langle (q_1-q_2)u,v\rangle_{L^2(\Omega_T)} \right|\leq C\|q_1-q_2\|_{L^{\infty}(0,T;L^p(\Omega))}\|u\|_{L^2(0,T;\widetilde{H}^s(\Omega))}\|v\|_{L^2(\Omega_T)}
        \end{split}
    \end{equation}
    for all $u\in L^2(0,T;\widetilde{H}^s(\Omega))$ and $v\in L^2(\Omega_T)$. 
    Next, by replacing $v_1$ by $v_1^k$ and $u_2-\varphi_2$ by $u_2^k-\varphi_2^k$ in \eqref{eq: identity for limit in uniqueness proof} we have
    \begin{equation}
    \label{eq: identity for each k}
         0=\int_{\Omega_T} (\partial_{\tau}f_1(v_0)-\partial_{\tau}f_2(v_0))v_1^k(u^k_2-\varphi^k_2)^{\star}\,dxdt
    \end{equation}
    for all $k\in\N$. Using \eqref{eq: estimate for difference}, we see that in the limit $k\to\infty$ the identity \eqref{eq: identity for each k} converges to
    \begin{equation}
    \label{eq: identity difference with v0}
         0=\int_{\Omega_T}  (\partial_{\tau}f_1(v_0)-\partial_{\tau}f_2(v_0))\Psi_1\Psi_2^\star\,dxdt.
    \end{equation}
    Next, we replace $v_0$ by $v_0^k$ to get
    \begin{equation}
         \label{eq: identity difference with v0}
         \begin{split}
         0&=\int_{\Omega_T} (\partial_{\tau}f_1(v^k_0)-\partial_{\tau}f_2(v^k_0))\Psi_1\Psi_2^\star\,dxdt\\
         &=\int_{\Omega_T} (\partial_{\tau}f_1(v^k_0-\psi_0^k)-\partial_{\tau}f_2(v^k_0-\psi_0^k))\Psi_1\Psi_2^\star\,dxdt.
         \end{split}
    \end{equation}
    Using Lemma~\ref{lemma: Continuity of Nemytskii operators}, we obtain that
    \begin{equation}
    \label{eq: continuity}
        \partial_\tau f_j\colon L^q(0,T;L^{r+2}(\Omega))\to L^{\frac{q}{r}}(0,T;L^{\frac{r+2}{r}}(\Omega))
    \end{equation}
    is continuous for any $r\leq q<\infty$. Next, recall that we have the embedding
    \begin{equation}
    \label{eq: embedding uniqueness proof}
        H^s(\R^n)\hookrightarrow L^{r+2}(\Omega)
    \end{equation}
    (see \eqref{eq: Hs to L r2} for the case $2s<n$). Therefore, as by assumption we have $0<r\leq 2$, we may from \eqref{eq: continuity} and \eqref{eq: embedding uniqueness proof} that $\partial_\tau f_j$ is continuous as a map from $L^2(0,T;H^s(\R^n))$ to $L^{\frac{2}{r}}(0,T;L^{\frac{r+2}{r}}(\Omega))$. Hence, by H\"older's inequality we can pass to the limit in \eqref{eq: identity difference with v0} and get
    \[
        \int_{\Omega_T} (\partial_{\tau}f_1(\Psi_0)-\partial_{\tau}f_2(\Psi_0))\Psi_1\Psi_2^\star\,dxdt=0.
    \]
    This implies that
    \[
        \partial_\tau f_1(x,\Psi_0(x,t))=\partial_\tau f_2(x,\Psi_0(x,t))\text{ in }\Omega_T
    \]
    for any $\Psi_0\in C_c^{\infty}(\Omega_T)$. This in turn allows us to conclude that
    \[
        \partial_\tau f_1(x,\rho)=\partial_\tau f_2(x,\rho)\text{ for all }(x,\rho)\in\Omega\times \R.
    \]
    Now, by the homogeneity of $f_j$ we can invoke Euler's homogeneous function theorem to conclude that \eqref{eq: equality nonlinearities} holds and we can finish the proof of Theorem~\ref{Thm: main nonlinear}.
    
    \end{proof}

	\medskip 
	
	\noindent\textbf{Acknowledgments.}
    The author gratefully acknowledges the support of the Swiss National Science Foundation (SNSF) under Grant No.~214500. The author thanks T.~Tyni and Y.-H.~Lin for helpful discussions. The author also thanks the handling editor and the anonymous reviewers for their valuable comments, which improved the presentation of the manuscript.

	\bibliography{refs} 
	
	\bibliographystyle{alpha}
	
\end{document}